\input amstex
\documentstyle{amsppt}
\magnification=\magstep0
\define\cc{\Bbb C}
\define\z{\Bbb Z}
\define\r{\Bbb R}

\define\N{\Bbb N}

\define\A{\Cal A}
\define\h{\Cal D}

\define\T{\Cal T}
\define\f{\Cal S}

\define\n{\Cal F}

\define\la{\lambda}

\define\e{\varepsilon}
\define\va{\varphi }
\define\CB#1{\Cal C_b(#1)}
\define\st{\subset }
\define\al{\alpha}
 \topmatter
 \title
Heat equation for weighted Banach space valued  function spaces
\endtitle
 \author
 Bolis Basit and Hans G\"{u}nzler
\endauthor
\abstract{ We study the homogeneous equation (*) $ u' = \Delta u$, $t > 0$, $u(0)=f\in wX$, where $wX$ is a weighted Banach space,  $w(x)= (1+||x||)^k$, $x\in \r^n$ with $k\ge 0$, $ \Delta$ is the Laplacian, $Y$ a complex Banach space and  $X$ one of the spaces  $   BUC (\r^n,Y)\} $, $ C_0 (\r^n,Y)$,
$ L^p (\r^n,Y)$, $1 \le p < \infty$. It is  shown that the mild solutions of (*) are still  given by  the classical Gauss-Poisson formula, a holomorphic $C_0$-semigroup}.
 \endabstract
\endtopmatter
\rightheadtext{Heat equations} \leftheadtext{ Basit and
 G\"unzler}
 \TagsOnRight

\document
\baselineskip=20pt

\head{\S 1.  Introduction, Notation and Preliminaries}\endhead

 In this note \footnote {AMS subject classification 2010:  Primary  {35K05, 47D06} Secondary {35K08,  46B99, 43A60}.
\newline\indent Key words: Laplace operator, holomorphic semigroups, weighted Banach space valued  function spaces.}  Example 3.7.6 of [1, p. 154] about solutions of the heat equation via holomorphic $C_0$-semigroups is extended to weighted function spaces and Banach space valued functions. Our treatment is different from [1, p. 154]: instead of using  Fourier transforms, direct methods are used.

\noindent   Let $w(x): =w_k (x)= (1+||x||)^k$ with $k\in \r_+= [0,\infty)$, $x=(x_1, \cdots,x_n)\in \r^n$, $||x||= (\sum_{k=1}^n x_k^2)^{1/2}$. Then $w\in C(\r^n)$ and

\

(1.1) \qquad $1 \le w(x+y) \le w(x)w(y) $, $w(y)\le w(x-y)w(x)$, $w(0)=1$,

\qquad \qquad\,\,\, $|w (x+y)/w(x)-1| \le w(y)(w(y)-1)$, \qquad    $x,y \in \r^n$.

\noindent Let $Y$ be a complex Banach space
and

\

(1.2) \qquad $w X= \{w g: g\in X\}$ with $X$ one of the spaces

\qquad \qquad\quad $   BUC (\r^n,Y) $, $ C_0 (\r^n,Y)$,
$ L^p (\r^n,Y)$, $1 \le p < \infty$.

 \noindent Then $w X $ is a Banach space  with norm $||f||_{w X}=||f/w||_{X}$ and a linear subset of $\f' (\r^n,Y)$, $wX$ is translation invariant, since X is and, with $f = wg$, $g \in X$, $f_h(x) : =
     f(x+h)$, one has $f_h/w = g_h w_h/w$ with $w_h/w  \in  BUC(\r^n,\r)$, using
     (1.1). For any  $f:\r^n \to Y$,   $|f| (x)= ||f(x)||$,
       $x\in \r^n$.

 \noindent     For $f\in w X$ and $\zeta \in \cc^+ :=\{ z\in \cc: \text {\,Re\,} \zeta > 0\}$ define (see Lemma 1.3)

\

(1.3) \qquad $(G(\zeta)f)(x):= (4\pi \zeta)^{-n/2} \int_{\r^n} f(x-y) e^{-||y||^2/4\zeta}\, dy$,  $x\in \r^n$.

\noindent Let
$\chi_{\zeta} (x)= (4\pi \zeta)^{-n/2}e^{-||x||^2/4\zeta}$, $\zeta \in \cc^+$,  $x\in \r^n$. Then $\chi_{\zeta}  \in \f (\r^n)$ if $\zeta \in \cc^+$,

\

(1.4) \qquad $(G(\zeta)f)= \chi_\zeta *f $, \,\, $\zeta \in \cc^+$,\,\,  $G(0)f=f$,\qquad $f\in w X$.

\noindent The function $\chi _{\zeta}$ is defined  and  $\chi'_{\zeta}=\frac {d \chi _{\zeta}}{d\zeta}$ exists  for  each $ \zeta \in \cc^+$, thus   holomorphic  on $\cc^+$. Moreover, $\chi^{(k)}_{\zeta}=\frac{d^k \chi _{\zeta}}{d\zeta^k} \in \f(\r^n)$ for each $\zeta\in \cc^+$, $k \in \N_0$.

\

(1.5) \qquad  $I= I (\zeta)= ((4\pi \zeta)^{- n/2}\int_{\r^n} e^{-(||x||^2/4\zeta)} \, dx= 1$ for each  $\zeta \in \cc^+$.

\noindent Indeed, $ I (\zeta)$ is holomorphic on $\cc^+$ with $I=1$ on $(0,\infty)$. It follows $I=1$ on $\cc^+$ by the identity theorem for complex valued holomorphic functions.

\noindent Also, for $\zeta = r e^{i\phi}$,   $0 \le |\phi| < \al  <  \pi/2$, $r > 0$, for any $x\in \r^n$

\

(1.6) \qquad $ |\chi_{\zeta} (x)|= (4\pi r)^{- n/2} e^{- (||x||^2 \cos \phi )/4r}  <  (4\pi r)^{- n/2} e^{- (||x||^2 \cos \al )/4r}$.

\

(1.7) \qquad Fourier transform $\widehat{\chi_{\zeta}}(x)=e^{-\zeta\, ||x||^2}$, $x\in \r^n$, $\zeta \in \cc^+$.

\noindent Indeed, it is enough to prove the case $n=1$. We have

$\widehat{\chi_{\zeta}}(y)=  e^{-\zeta y^2} I(\zeta,y)$, where

$   I(\zeta,y)= (4\pi \zeta)^{-1/2} \int_{\r}  e^{-(x+2i\zeta y)^2/4\zeta} \,dx$, $y\in \r$, $\zeta \in \cc^+$.

\noindent With
$F(x,y) : = e^{-(x+2i \zeta y)^2/4\zeta}$,

 $\frac{\partial}{\partial y} \int_{\r} F(x,y)\, dx  =$
$\int_{\r}\frac{\partial}{\partial y}   F(x,y)\, dx $
$= \int_{\r} 2i \zeta \frac{\partial}{\partial x} F(x,y) dx =$

\qquad \qquad \qquad\qquad $= 2i\zeta \lim_{N\to \infty}
(F(N,y) - F(-N,y)) = 0$,

\noindent  so $I(\zeta,y) = I(\zeta,0)$, $= 1$ for  $\zeta$
real  $> 0$  (e.g. [3,  p. 274, Beispiel 1]), then for  $\zeta \in \cc^+$
since $I$ is holomorphic there.

\proclaim{Lemma 1.1} If $f\in w X$ respectively $wf \in L^p (\r^n,\cc)$ with $1 \le p <\infty$, then $||(f_y-f)/w||_{X}\to 0$  respectively $||w(f_y-f)||_{L^p}\to 0$ as $ y \to 0$.
\endproclaim
\demo{Proof} Let $f=w g$, where $g\in X$. Then $||f_y-f||_{wX}=|| w_y g_y -wg||_{wX}= || (w_y-w) g_y+ wg_y -wg||_{w X}\le || (w_y-w) g_y||_{wX} + || wg_y -wg||_{wX} = || (w_y/w -1) g_y||_{X} + || g_y -g||_{X}\to 0$ as $y\to 0$ with (1.1) since   $|| g_y -g||_{X}\to 0$ as $y\to 0$. The second case follows similarly.\qquad\qquad  $\square$
\enddemo

\proclaim{Lemma 1.2}
 (A)  If  $f\in w L^p(\r^n,Y)$, $wg\in L^q(\r^n,\cc)$  with $1/p + 1/q = 1$ and
    $1\le p\le \infty $, then
     $(g*f)(x)$ exists as a Bochner integral for all $x\in \r^n$, and
          $g*f\in wBUC(\r^n,Y)$;
    if additionally  $ 1<p< \infty$ or $f\in w C_0(\r^n,Y)$ and $q=1$, then $g * f  \in w C_0 (\r^n,Y)$.

 (B)  If  $f\in wL^p(\r^n,Y)$,  $wg \in L^1(\r^n,\cc)$ with $1\le p\le \infty$,
         then   $g*f(x)$ exists as Bochner integral almost everywhere in $\r^n$ and
               $  g*f  \in  w L^p(\r^n,Y) $.
\endproclaim
\demo {Proof} (A) Since

\

(1.8)\qquad  $||f(y)g(x-y)|| = ||(f/w)(y)|| |(wg)(x-y)|  (w(y)/w(x-y))  \le
     $

  \qquad \qquad    $|f/w|(y)$
      $|wg|(x-y) w(x)$,

 \noindent  (1.1),  $|f/w| \in  L^p (\r^n)$ and $|wg|(x-\cdot) \in L^q (\r^n)$,
 with the H\"{o}lder inequality [5, p. 34, Proposition 2] one has $f(\cdot)g(x-\cdot)
 \in L^1(\r^n)$ for all $x \in \r^n$,

 \
 
  (1.9) \qquad $||g * f (x)|| \le w(x)||wg||_{L^q}|| f/w||_{L^p} $.

  \noindent  With this

   $||g * f (x+y)-g * f (x)||\le  $
   $w(x) ||f/w||_{L^p} ||w(g _y- g)||_{L^q}$,

    $||g * f (x+y)-g * f (x)||\le   $
   $w(x) ||(f_y-f)/w||_{L^p} ||wg||_{L^q}$.

\noindent By Lemma 1.1,  $||w(g _y- g)||_{L^q}\to 0$ respectively $||(f_y-f)/w||_{L^p}\to 0$ as $y\to 0$ if  $ 1 \le q < \infty$ respectively $ 1 \le p < \infty$. It follows
   $g * f\in w BUC (\r^n,Y)$  if  $ 1 \le p,  q \le \infty$.

\noindent If  $ p > 1, q < \infty$ or $f\in w C_0(\r^n,Y)$ and $q=1$, then  $(w|g|)* (|f|/w)\in C_0(\r^n)$ by [1,Proposition 1.3.2 b), d),  p. 22 ]. It follows  $g * f  \in w C_0 (\r^n,Y)$.

(B) By Young's inequality [5, p. 29], $(w|g|)* (|f|/w) \in L^p(\r^n)$. So, $(w|g|)* (|f|/w) (x))$ is finite almost everywhere on $\r^n$. This,  measurability of $ g(x-\cdot)  f(\cdot)$  and (1.8) imply $g * f (x)$  exists as a Bochner integral almost everywhere on $\r^n$. The above $(w |g|)*(|f|/w)\in L^p (\r^n,\cc)$  and (1.8) give  $g * f\in w L^p (\r^n,Y)$.\qquad\qquad  $\square$
\enddemo

\proclaim{Lemma 1.3} Let  $f\in w X$, $G(\zeta)$ defined by (1.3n) and  $g =\chi_{\zeta}$ or $\chi'_{\zeta}:= \frac{d \chi_{\zeta}}{d\zeta}$, $\zeta\in \cc^+$.

(i) $g * f (x)$  exist as a Bochner integral  for all $x\in \r^n$ and $g * f \in w BUC(\r^n,Y) \cap w X$;  if additionally  $ 1<p< \infty$ or $f\in w C_0(\r^n,Y)$,  then  $ g * f \in w C_0 (\r^n,Y)\cap w X$.

(ii) $G(\zeta)\in L(wX)$.

 (iii) If   $0 < \al <\pi/2$, then

\

(1.10) \qquad $\lim _{0 \not = \zeta \to 0, |\text {arg\,} \zeta| <\al} ||\chi_{\zeta} * f-f||_{w X}= 0$.
\endproclaim
\demo{Proof}  (i) Since  $wg \in L^q (\r^n,\cc)$ for each $1 \le q \le \infty$, (i) follows by Lemma 1.2.

(ii) The operator $G(\zeta): wX\to wX$  defined by $G(\zeta)f :=\chi_{\zeta}*f$ is linear  and bounded by (1.9).

 (iii) With $y= |\zeta|^{1/2}z$ and $\theta= \frac{\zeta}{|\zeta|}$, it follows by (1.5)

$\chi_{\zeta} * f (x)-f(x)= \int_{\r^n} [f(x-y)- f(x)] \chi_{\zeta} (y)\, dy= \int_{\r^n} [f(x-|\zeta|^{1/2}z)- f(x)] \chi_{\theta} (z)\, dz$.

\noindent Case $X= BUC (\r^n,Y)$,  $ C_0(\r^n,Y)$. Let $\e > 0$. Since $w  \chi_{\theta} \in L^1 (\r^n)$, then  using  (1.1),  for  $0 < |\zeta|  \le 1$, $|\text{arg\,\,}\zeta| < \al$  there exists $c= c(\e,\al) > 0$ independent of $\zeta$,  such that

$I_1=$  sup $_{x\in \r^n}\frac{1}{w(x)}  \int_{||z|| \ge c} ||f(x- |\zeta|^{1/2}z)- f(x)|| |\chi_{\theta} (z)|\, dz \le 2 ||f||_{wX}\times$

 $\int_{||z|| \ge c} w (z) |\chi_{\theta} (z)|\, dz  \le    2 ||f||_{wX} (4\pi)^{-n/2} \int_{||z||>c} w(z)e^{-||z||^2(\cos\, \al)/4} dz
             <   \e$.

\noindent Then for  the above  $\zeta$

\noindent $ ||\chi_{\zeta} * f (x)-f(x)||_{w X}\le $
 $I_1 +$   sup $_{x\in \r^n} \frac{1}{w (x)} \int_{||z|| \le c} ||f(x- |\zeta|^{1/2}z)- f(x)|||\chi_{\theta} (z)|\, dz \le $

\noindent  $ I_1 + $ sup $_{x\in \r^n, ||z||\le c} \,\,\, \frac{||f(x- |\zeta|^{1/2}z)- f(x)||}{w(x)} (4\pi)^{-n/2} \int_{\r^n} e^{-||z||^2 (\cos \al)/4 }\, dz= I_1 +I_2$.

\noindent   Using Lemma 1.1, there is $\delta >0$ such that $I_2 \le \e$ if $|\zeta|^{1/2} c < \delta$. It follows $I_1+I_2 \le   2 \e$  if  $ 0<|\zeta|^{1/2}< \delta/c  $ and $|\text{arg\,} \zeta| \le \al$.

 Case $X= L^p$: By  (i)
 $\chi_{\zeta}*f  \in   w L^p(\r^n,Y) \cap BUC(\r^n,Y)$.
 For $  \zeta \in \cc^+$ with  $y= |\zeta|^{1/2}z$  using  the Minkowski inequality [3,  p. 251, A 92]

\noindent $||\chi_{\zeta} * f -f||_{{w} L^p}= [\int_{\r^n}\frac{||\int_{\r^n}  [f (x-y)-f(x)] \chi_{\zeta} (y)\, dy ||^p}{w^p_k (x)} dx]^{\frac{1}{p}}\le  $

\noindent $ \int_{\r^n}[\int_{\r^n} \frac{  ||f (x-y)-f(x)||^p} {w^p_k (x)}\, dx]^{\frac{1}{p}} |\chi_{\zeta} (y)|\, dy$
 $ = \int_{\r^n}[\int_{\r^n} \frac{ ||f (x-|\zeta|^{1/2}z)-f(x)||^p} {w^p_k (x)}\, dx]^{\frac{1}{p}} |\chi_{\theta} (z)|\, dz$.

  \noindent By Lemma 1.1,  $\int_{\r^n} \frac{ ||f (x-|\zeta|^{1/2}z)-f(x)||^p} {w^p_k (x)}\, dx]^{\frac{1}{p}}\to 0$ as $|\zeta|\to 0$  for each $z\in \r^n$.  So,  by the dominated convergence theorem as in Lemma 1.3.3 (b) of [1, p. 23] we get the statement since

   \noindent $|\chi_{\theta} (z)| < (4 \pi)^{-n/2} e^{-||z||^2 (\cos \al)/4}= : F(z)$  by (1.6), $||f_{-|\zeta|^{1/2}z}||_{wX}  \le  w(z) ||f||_{wX}$  and $wF\in L^1 (\r^n)$,
         if $z \in \r^n$, $|arg \zeta | < \al $, $0 < |\zeta| \le 1$.\qquad\qquad  $\square$
\enddemo

\head{\S 2.  Main results}\endhead

\proclaim {Theorem 2.1} For $w X$  of (1.2), the $G$  of (1.3) is a  holomorphic $C_0$-semigroup of angle $\pi/2$ on $w X$. Its generator is the Laplacian $\Delta_{wX} :=\Delta$ on $w X$ with  domain:

$D(\Delta_{wX})= \{f\in wX:\text {\,distribution-} \Delta f\in w X\}$, $\Delta = \sum_{j=1}^n \partial ^ 2  /\partial x^2_j$

\noindent where we identify $wX$ with  a subspace of $\f' (\r^n,Y)$.
\endproclaim

\demo{Proof} (a): We have  $\chi_{\zeta} \in \f (\r^n)$  for $\zeta\in \cc^+$ and

$\frac{d\chi_{\zeta}}{d\zeta} (x) =\Delta \chi_{\zeta} (x)$ for $\zeta \in \cc^+$, $x\in \r^n$.

\noindent Moreover, by Lemma 1.3  $G(\zeta)f= \chi_{\zeta}*f \in w X$,  $||\chi_{\zeta}*f- f||_{w X} \to 0$ as in  (1.10) for all $\zeta\in \cc^+$, $f\in wX$ and $G(\zeta)\in L(wX)$. Then  $\widehat{G(\zeta)f}= \widehat{\chi_{\zeta}}\cdot \widehat{f}$ follows as in [1, p. 154]. By (1.7),   $\widehat{\chi_{\zeta_1+\zeta_2}}=\widehat{\chi_{\zeta_1}}\widehat{\chi_{\zeta_2}}$. So, $G(\zeta_1+\zeta_2)=G(\zeta_1)G(\zeta_2)$, $\zeta_1, \zeta_2\in \cc^+$.
  This means that $G$ is a $C_0$-semigroup on $wX$.

(b)  Holomorphy of $G : \cc^+ \to  L(wX)$.   By [1, Proposition A.3, (ii) $\Rightarrow$ (i)],
  it is enough to show that for any  $f \in wX$  with  $U(\zeta)= G(\zeta)f$ the $U$ is holomorphic
  on $\cc^+ $.
Now, again   by [1, Proposition A.3], holomorphy of the function $\zeta \to w \chi_{\zeta}$  defined
on $\cc^+ $ with values in $L^1(\r^n)$ follows, since  the complex  valued
$ F(\zeta)=\int_{\r^n} w (x)\chi_{\zeta} (x) g (x)\, dx$  is continuous  for each $g \in L^{\infty} (\r^n)$ and   by Morera's
theorem [4, p.75],   Fubini  and (1.6) it is holomorphic.
 So to fixed z there exists $\psi$ in $ L^1 (\r^n)$ with  $w(\frac{\Delta \chi_{\zeta}}{\Delta \zeta}) \to \psi$  in  $L^1 (\r^n)$; so
there are $\zeta_n \to  \zeta$ with  $\frac {\chi_{\zeta_n} - \chi_{\zeta}}{\zeta_n -\zeta}
\to  \psi/w$  almost everywhere   on $\r^n$; with the holomorphy of $\chi_{\zeta}(x)$ for each $x\in  \r^n $   one gets $\psi/w =
\chi'_{\zeta}$ almost everywhere and

\

 (2.1)\qquad  $||(\frac{\Delta \chi_{\zeta}}{\Delta \zeta}- \chi'_{\zeta}) w||_{L^1}= \int_{\r^n}|\frac{\Delta \chi_{\zeta} (x)}{\Delta \zeta}- \chi'_{\zeta} (x)| w(x)\, dx \to 0$ as $0 \not =\Delta \zeta \to 0$.

\noindent Since  $w \chi'_{\zeta}  \in   L^q (\r^n)$ for all $ q\ge 1$,  $\chi'_{\zeta}  * f (x)$ exists with H\"{o}lder's inequality
as  a Bochner integral for  all $x \in \r^n $.
By Lemma 1.3, $\frac{\Delta U(\zeta)}{\Delta \zeta},\chi'_{\zeta}* f\in wX$, $\zeta, \zeta+\Delta \zeta\in \cc^+$, $\Delta \zeta \not = 0$. We have

  $\frac{\Delta U(\zeta)}{\Delta \zeta}- \chi'_{\zeta}* f= (\frac{\Delta \chi_{\zeta}}{\Delta \zeta}- \chi'_{\zeta})*f$ and  using Young's inequality [5, p. 29]

$||(\frac{\Delta \chi_{\zeta}}{\Delta \zeta}- \chi'_{\zeta})*f||_{wX}= ||(1/w)\int_{\r^n}(\frac{\Delta \chi_{\zeta}}{\Delta \zeta}- \chi'_{\zeta})(\cdot-y)f(y)\, dy||_{X} \le$

$ ||\int_{\r^n}|(\frac{\Delta \chi_{\zeta}}{\Delta \zeta}- \chi'_{\zeta})(\cdot-y)| w(\cdot-y)(||f(y)||/w(y))\, dy||_{X}\le  $

$||(\frac{\Delta \chi_{\zeta}}{\Delta \zeta}- \chi'_{\zeta}) w||_{L^1}||f/w||_X$.

\noindent   With (2.1), holomorphy of $U$ on $\cc^+$ follows, and

\

 (2.2) \qquad   $G'(\zeta)f = (G(\zeta)f)' = \chi_{\zeta}' * f$,   $\zeta \in  \cc^+$,  $f \in  wX $.

 (c) Let $f$,  distribution $\Delta f  \in w X$. We have $\frac{\partial \chi_t} {\partial t}=\Delta_x \chi_t$ on  $(0, \infty)\times \r^n$, $\Delta_x = \sum_{j=1}^n (\frac{\partial} {\partial x_j})^2$. So by (2.2), in $\f' (\r^n,Y)$, $t > 0$,

\

 (2.3) \qquad   $\frac{d G(t)f}{dt}= \frac{d (\chi_t *f)}{d t}= \frac{d \chi_t} {d t}*f=$

 \qquad \qquad \quad $ (\Delta \chi_t )* f = \Delta G(t)f =\chi_t *(\Delta f)= G(t) \Delta f$.

\noindent Let $A|D(A)$ be the generator of  the $C_0$- semigroup  $G : \r_+ \to L(wX)$,
defined by  Proposition 3.1.9 g) of [1,  p. 115 ]; let $\Delta$ be the Laplace operator applied
to  $S  \st  \h' : = \h'(\r^n,Y) $;
 with  $wX \st \h'(\r^n,Y)$,

$D : =\{f \in wX : \Delta_{wX} f  \in wX\}$ and $ \Delta_{wX}: = \Delta | D$
are well defined.  We show

\

(2.4)\qquad   $D(A) = D$,  $A = \Delta_{wX}$.

(c.1)\qquad $D \st  D(A)$,   $A = \Delta_{wX}$  on  $D$:

\noindent  If  $f \in D$,
$G(\cdot)f \in C([0,\infty),wX)$  by (a), with (2.3) and   $g : =  \Delta_{wX} f  \in wX$  one has
 $G(f)f  - f =  \int_0^t  (\frac{d}{ds})(G(s)f)\, ds = \int_0^t  G(s)g \, ds$, $t \in \r_+ $.
 With Proposition 3.1.9 f) of [1,  p. 115 ]  one gets $f \in D(A)$  and  $Af = g = \Delta_{wX} f $.

(c.2) \qquad   $D(A) \st D$ :

\noindent  With  $F(t) : = (1/t) \int_0^t  G(s)f\, ds$, $t> 0$, $f \in D(A)$,
    one has  $F(t)\to f$  in $wX$ as  $t\to 0$, since  $G(t)f\to f$ in $wX$  by (a).
    $F(t)\to f$ in $wX$ implies $F(t)\to f$  in $L^1_{loc}(\r^n,Y)$, so
$(\Delta F(t))(\va) =  F(t)(\Delta_{x} \va)\to f (\Delta_{x} \va) = (\Delta f)(\va)$
for  $\va \in \h(\r^n,\cc)$.

Now by (2.5) below one  has  $\Delta F (t) = (1/t)(G(t)f -f) $;
by definition of $D(A)$ and  Proposition 3.1.9 g) [A., p. 115], $(1/t)(G(t)f -f)\to$ some $g$
in $wX$, so in $L^1_{loc} (\r^n,Y)$, so $(1/t)(G(t)f-f)(\va)\to f(\va)$.
together one gets $\Delta f = g $, $\in wX$, that is  $f \in D$.
With  (c.1) this gives  (2.4).
It remains to show

\

(2.5) \qquad   $\Delta \int_0^t G(s)f ds  =  G(t)f - f$,  $f  \in  wX $.

\noindent For this, with $f\,\in\, wX$, with Lemma 1.3 define $\beta(t,x) : =
(\chi_t*f)(x)$, $(t,x)\,\in\, M : = (0,\infty) \times \r^n$.
With Lebesgue's Dominated Convergence theorem and analogs of (1.6) for the
derivatives of $ \chi_t$ one gets inductively $\beta\,\in\, C^{\infty}(M,Y)$, with

\

(2.6) \qquad $ \partial \beta/\partial t = (\chi'_t)*f = (\Delta_x \chi_t)*f = \Delta_x\beta$.

\noindent If  0 $< \e < t$,  $\Psi _{\e} (t,x) : = \int^t _{\e}\beta(s,x) ds, x \in \r^n$, is
well defined with  $\Psi _{\e}\,\in\, C((\e,\infty)\times \r^n,Y)$, $\Psi _{\e}(t): = \Psi _{\e}(t,\cdot)$
$\in C(\r^n,Y) \st  \h' (\r^n,Y)$  if $t  > \e$. If $\va\in \h (\r^n)$, all the following integrals
exist (even as Riemann integrals), with twice Fubini, partial
integration and (2.6) one has

$(\Delta \Psi _{\e}(t))(\va) = \int_{\r^n} \Psi _{\e}(t, x)(\Delta_x \va)(x)\,dx =
 \int_{\r^n}\int^t _{\e} \beta(s,x)(\Delta_x \va)(x)ds\, dx =$

$ \int^t _{\e} \int_{\r^n}
\Delta_x\beta (s,x) \va(x)dx ds = \int^t _{\e} \int_{\r^n} (\partial/\partial s)\beta(s,x)\va(x)
dx \, ds =  $

$\int_{\r^n} (\int^t _{\e} (\partial/\partial s)\beta(s,x)ds)\va(x) dx =  $
$\int_{\r^n}
(\beta(t,x) -\beta(\e,x))\va(x) dx $.

\noindent This implies

\

(2.7)  \qquad  $ \Delta \Psi _{\e}  =  G(t)f  -  G(\e)f$, $\in\,  wX $.

\noindent $G(\cdot)f :\r_+\,\to\,wX$ is continuous, so $\int^ t _{\e} G(s)f\, ds\,\to\,
\int^t_0 G(s)f\, ds $  as  $\e\,\to\, 0  $.
Furthermore, the Riemann sums $\Sigma_m : = \sum^m_1 (G(s_j)f)\, (s_j-s_{j-1})
\to \int^t _{\e} G(s)f \,ds$   in $wX$  as  $m\to \infty$, $s_j = \e + j(t-\e)/m$.
Similarly  $\Sigma_m(x): = \sum^m_1 \beta(s_j,x)(s_j - s_{j-1})\,\to\,
\int^t _{\e}\beta (s,x)\,ds = \Psi _{\e}(t,x)$  in $Y$  as $ m \to \infty$, for each $x\,\in\, \r^n$.

\noindent If $K$ is compact $\st \r^n$, then sup $\{||\Sigma_m(x)||: m\,\in\N, x\,\in\, K\} < \infty$,
so

$\int_{\r^n} \Sigma_m(x)\va(x)\, dx\,\to\, \int_{\r^n} \Psi _{\e}(t,x) \va (x)\, dx$ for $\va\,\in\, \h (\r^n)$.

\noindent As above, $\int^t _{\e} G(s)f\,ds = \Psi _{\e}(t)$  follows, and then
$\Psi _{\e}(t)\to \int^t_0 G(s)f\, ds$  in $wX$. Therefore
$(\Delta \Psi _{\e}(t))(\va) = \int_{\r^n} \Psi _{\e}(t) \Delta_x \va\,dx \,\to\,
\int_{\r^n}(\int^t_0 G(s)f \, ds)\Delta_x \va\, dx=$

\noindent $  (\Delta \int^t_0 G(s)f\, ds)
(\va) $ as  $\e\,\to\, 0$;
  since  $G(\e)f\,\to\, f$, (2.7) implies (2.5).\qquad\qquad  $\square$
\enddemo

\proclaim {Corollary 2.2} All mild solutions $u : [0,\infty) \to wX$  of

(a) $ u' = \Delta_{wX} u$
   are given by $u(t) = G(t)f$  with $f \in  wX$, they are $C^1$-solutions on
   $(0,\infty)$, $\in  C^{\infty}(M,Y)$ and classical solutions of

   (b)  $\frac{\partial{u(t,x_1, \cdots,x_n)}} {\partial t}= \sum^n_1 \frac{\partial^2{u(t,x_1, \cdots,x_n)}} {\partial x^2_j}$  on $M : = (0,\infty)\times \r^n$  with

    (c)  $u(t,\cdot) \to  f$ in  $wX$ as $t \to 0 $.

\noindent     Conversely, any classical solution of (b) with (c) defines a mild solution
    of (a) on $[0,\infty)$.
    \endproclaim
For the proofs of most of  this  see [1, Corollary 3.7.21].

\noindent {\bf Remark 2.3}.  Since the Gauss-Poisson formula  (1.4) for $G$ defines by Theorem 2.1
a holomorphic $C_0$-semigroup with generator $A = \Delta_{wX}$, with Corollary 2.2
the results of  [2, Theorems 5.2/6.3, Examples 6.2] can be applied to the heat
equation.

\indent School of Math. Sci., P.O. Box No. 28M, Monash University,
 Vic. 3800.

\indent E-mail "bolis.basit\@monash.edu".

\indent Math. Seminar der  Univ. Kiel, Ludewig-Meyn-Str., 24098
Kiel, Deutschland.

\indent E-mail "guenzler\@math.uni-kiel.de".

\enddocument

Lemma. Let $f\in w L^p (\r^n,Y)$, $g\in \f(\r^n)$ and $F(y) =f_{-y} g(y)$, $y \in \r^n$. Then $F$ is Bochner integrable on $\r^n$ and

 $\int_{\r^n} F(y) = f*g \in w L^p (\r^n,Y)\cap wBUC(\r^n,  L^p (\r^n,Y))$.

\

(2.5) \qquad   $\Delta \int_0^t G(s)f ds  =  G(t)f - f$,  $f  \in  wX $

Proof. By Lemma 1.3, $ G'(s)f,  G(s)f\in BUC(\r^n,Y)\cap wX$ for each  $ s > 0$. It follows $G'(\cdot)f \in C((0,\infty],wX)$;  so Bochner integrable on $[\e,t]$, $0 < \e < t < \infty$ and
$ \int_{\e}^t G'(s)f \, ds= G(t)f- G(\e)f  $. By (1.12), $\int_{\e}^t G'(s)f \, ds = \int_{\e}^t \Delta G(s)f) \, ds $. Since
$\Delta $ is closed,   [A, Prop. 1.1.7, p. 11] gives  $ \Delta \int_{\tau}^t  (G(s)f) \, ds = \int_{\e}^t (\Delta (G(s)f)) \, ds  =  G(t)f - G(\e)f$,  $f  \in  wX $. Taking $0 < \tau \to 0$, (2.5) follows by (1.10).

\

(2.5) \qquad   $\Delta_{wX} \int_0^t G(s)f ds  =  G(t)f - f$,  $f  \in  wX $

Proof. By Lemma 1.3, $ G'(s)f,  G(s)f\in BUC(\r^n,Y)\cap wX$ for each  $ s > 0$. It follows $ (\int_{\tau}^t G(s)f \, ds) (x) = \int_{\tau}^t (G(s)f)(x) \, ds$ and  $ (\int_{\tau}^t G'(s)f \, ds) (x) = \int_{\tau}^t (G'(s)f)(x) \, ds$

\noindent $ = (G(t)f) (x)-(G(\tau)f) (x)$.  By (1.12),  $\Delta_{wX}$ is closed and [A, Prop. 1.1.7, p. 11], $\int_{\tau}^t (G'(s)f)(x) \, ds = \int_{\tau}^t (\Delta_{wX} (G(s)f))(x) \, ds =\Delta_{wX} \int_{\tau}^t  (G(s)f)(x) \, ds$. This gives $\Delta_{wX} \int_{\tau}^t G(s)f ds  =  G(t)f - G(\tau)f$,  $f  \in  wX $. Taking $0 < \tau \to 0$, (2.5) follows by (1.10).

\ref\no5\by H. Reiter, Classical Harmonic Analysis and Locally Compact Groups, Oxford Math. Monographs, Oxford Univ.,  1968
\endref

 \noindent  For $\psi \in \f(\r^n)$,  $\int_0^t \frac{d\chi_s}{ds} * \psi \, ds=\chi_t *\psi- \psi $.

\noindent Now, we show that

(2.4) \qquad  $\Delta \int _0^tG(s)f\, ds = G(t)f-f$ for $f \in w X$ and $t >0$.

\noindent Then by [A, Corollary 3.1.13] and the above,  the generator of $G$ is the Laplacian with domain $D(\Delta_{wX})$.

Let $\psi\in \f(\r^n)$. Then Fubini's theorem gives

$\langle \Delta \psi, G(s)f \rangle =\langle G(s) \Delta \psi,f \rangle $ and

$\langle \psi,\Delta \int_0^t G(s)f \, ds\rangle = \langle \Delta \psi,\int_0^t G(s)f \, ds\rangle =\int_0^t \langle \Delta \psi, G(s)f \rangle \,ds$

$ \int_0^t \langle G(s)\Delta \psi, f \rangle \,ds= \langle\int _0^t G(s)\Delta \psi\, ds, f \rangle = \langle G(t) \psi-\psi,f \rangle= \langle \psi, G(t) f-f \rangle$.

\noindent This proves (1.12).

 With   $f = wg  \in  w L^p$,  $g \in L^p$,  and $\zeta\in \cc^+$ fixed, consider

  $F(x,y) : = \chi_{\zeta}(x)w(x) $
   $h(x,y) g(y-x)$, where

     $h(x,y) : = w(x-y)/(w(x)w(y))$,   $0< h(x,y) \le 1$ for all  $x,y \in \r^n$.

\noindent    Since
$\chi_{\zeta}, w, h$ are continuous on $\r^n$  respectively  $\r^n\times \r^n$,
one checks easily that all assumptions for the continuous Minkowski inequality [3, Int. p.251, A 92]"
         are fulfilled, so one
gets :
 $(\chi_{\zeta}*f)(y)$ exists as a Bochner integral for almost all $y \in \r^n$,

 $||\chi_{\zeta}*f||_{wL^p} \le  ||f||_{w L^p} \int_{\r^n}w(x)
  \chi_{\zeta}(x) dx$, so $G(\zeta)f : = \chi_{\zeta}*f$ defines a $G(\zeta) \in L(w L^p(\r^n,Y))$.

 If $wg$ is bounded and $wg\in L^1$, then  $wg \in L^q (\r^n,\cc)$ for all $1\le q \le \infty $,
            so  (A) applies.
            This is the case if   $g = \chi_{\zeta} $  or  $\chi'_{\zeta}$, $\zeta$ fixed  $\in \cc^+ $.

We have

(1.10)\qquad $\chi_{\zeta} * f (x) = \int_{\r^n}  f (y) \chi_{\zeta} (x-y)\, dy = \int_{\r^n} \frac {f (y)}{w(y)} w (y)\chi_{\zeta} (x-y)\, dy$.

(a) $X=  BUC (\r^n,Y)$ (the case $X=  C_0 (\r^n,Y)$ follows similarly):

$||\chi_{\zeta} * f (x)|| \le ||f/w ||_{\infty}\int_{\r^n}  w (y)|\chi_{\zeta} (x-y)|\, dy\le w(x)||f ||_{wX}||w |\chi_{\zeta}| ||_{L^1}$.

This implies,  $||\chi_{\zeta} * f (x+y)-\chi_{\zeta} * f (x)||\le w(x)||(f_y-f)/w ||_{\infty}||w |\chi_{\zeta}| ||_{L^1}$.

\noindent  Thus  with  $F :=\chi_{\zeta} * f$ the $F/w$ exists everywhere,  is bounded and  $ (F/w)_h -F/w = (F_h-F)/w + F_h (w/w_h-1)w \in  BUC(\r^n,Y)$ with Lemma 1.1 and (1.1).

 Let  $f\in w X$ and $wg \in L^1 (\r^n)$.

(i) If $X=L^{\infty}(\r^n,X)$, then  $g * f (x)$  exists as a Bochner integral,  $||g * f (x)||\le w(x) ((w|g|)* (|f|/w) (x)) $ for all $x\in \r^n$ and  $g * f  \in w BUC(\r^n,Y)$. Moreover, if $X=C_0(\r^n,Y)$, then $g * f  \in w C_0 (\r^n,Y)$.

(ii) If  $X=L^p(\r^n,Y)$, $1\le p <\infty$,  then $g * f (x)$  exists as a Bochner integral,  $||g * f (x)||\le w(x) ((w|g|)* (|f|/w) (x)) $ almost everywhere on   $ \r^n$, $g * f  \in w X$.

We have

\qquad $g * f (x) = \int_{\r^n}  f (y) g (x-y)\, dy = \int_{\r^n} \frac {f (y)}{w(y)} w (y)g (x-y)\, dy$.

\noindent  Using (1.1),

(1.*) \qquad $||g * f (x)|| \le  \int_{\r^n} \frac {||f (y)||}{w(y)} w (y) |g (x-y)|\, dy \le $

\qquad \qquad\,\, \, $ w(x) \int_{\r^n} \frac {||f (y)||}{w(y)} w (x-y) |g (x-y)|\, dy= w(x)((w|g|)* (|f|/w) (x)) $.

\noindent (i)  By (1.*),

$||g * f (x)||
  \le $
 $w(x) ||\frac {|f |}{w}||_{X} \int_{\r^n}  w (x-y)|g (x-y)|\, dy = w(x) ||f/w||_X ||wg||_{L^1}$.

  \noindent This  and    $ y^*\circ ((g * f )/w) \in BUC(\r^n)$ for each $y^*\in Y^*$ ([Reiter, (6), p. 85])  imply that $g * f (x)$ exists as a Bochner integral  for all $x\in \r^n$.
  As in (1.*),

   $||g * f (x+y)-g * f (x)||\le w(x) ((|f|/w)*(w|g _y- g|))(x)) \le  $

   $w(x) |||f|/w||_{L^p} ||w|g _y- g|||_{L^q}$.

\noindent Since $||w|g _y- g|||_{L^q}\to 0$ as $y\to 0$ if  $ 1 \le q < \infty$, we conclude   $g * f\in w BUC (\r^n,Y)$  if  $ 1 \le q < \infty$.  The case $q=\infty$ can be proved similarly. If if $X=C_0(\r^n,Y)$, then $g * f  \in w C_0 (\r^n,Y)$ by (1.*) and $(w|g|)* (|f|/w)\in C_0(\r^n)$ by [A, p.22 b)].

(B) By Young's inequality [5, p. 29], $(w|g|)* (|f|/w) \in L^p(\r^n)$. So, $(w|g|)* (|f|/w) (x))$ is finite almost everywhere on $\r^n$. This,  measurability of $ y^*\circ ((g(x-\cdot)  f(\cdot) )$ for each $y^*\in Y^*$  and integrability of $|g(x-\cdot)  f(\cdot)|$ imply $g * f (x)$  exists as a Bochner integral almost everywhere on $\r^n$. Again by (1.A*), $g * f\in w L^p (\r^n,Y)$

\proclaim{Theorem} For $w X$, $G$ as above, $G$ is a  holomorphic $C_0$-semigroup of angle $\pi/2$ on $w X$.
\endproclaim
\demo{Proof} We show that the complex derivative of
$U(z) : = G(z)f  = \chi_z * f , U : \cc^+ \to wX$,  exists for all $z \in \cc^+$.
 By [A, Prop. A.3, (ii) $\Rightarrow$ (i)],  for $ G : \cc^+ \to  L(wX)$ to be holomorphic
  it is enough to show that for any  $f \in wX$ the above $U$ is holomorphic
  on $\cc^+ $.
 So in the following, $f$ is fixed  $\in\,\,wX$, $z$ is fixed (but also arbitrary) $\in\,\,\cc^+$,
  one has to show that  complex derivative $U'(z)$ exists, i.e.

(1)\qquad  $\lim (G(\eta)f - G(z))/(\eta -z)$ exists in $wX$  as $\eta  \to z$, $z\not = \eta \in\,\,\cc^+$.

\noindent  (1) follows if  the above limit   exists for all sequences $(\,\eta_n) \in \cc^+$ with
    $\eta _n \to  z$, and since $wX$ is complete, it is enough to show that

 (2) \qquad   $\lim ( H(\eta_n) - H(\eta _m)) = 0$  in $wX$ if $n, m \to \,\infty$, with

   \qquad \qquad $H(\eta) : = (G(\eta)f  - G(z)f)/(\eta -z)$;

   \noindent  this in turn follows if

   (3)\qquad  $\lim (H(\eta _n) - H(\xi_n)) = 0 $ if $\eta _n \to z$ and
    $\xi_n \to z$.

\noindent  By definition of $G$, one has $H(\eta) = ((\chi_{\eta} -\chi_z)/(\eta -z))* f$,
    which by Lemma 1.2 is well defined and  $\in\,\,wX $ (always $\eta \, \not = z$, in
   $\cc^+$);
   Now $\chi_{\zeta }(x)$ has for fixed $x \in\,\,\r^n$  as a function  $\cc^+ \to \cc$ has everywhere in $\cc^+$
a complex derivative which I denote by $\chi_{\zeta}'(x)$,
 so  in  $ H(\eta _n) - H(\xi_n)  =  v_n  * f$  the   $(v_n)(x) \to 0$ as  $n\to\,\infty$,
 for every $x  \in\,\,\r^n$.

  This will be used to show  (3), with the aid of the dominated convergence theorem, so one has to get a suitable dominating
   function.
To the fixed $z \in\,\,\cc^+$  there exist positive real $a$ and $\al \,$ with $ a<1$ and
$\al \, < \pi/2$, so that the disk  $D : = \{ \eta  \in\cc : |\eta  - z| \le  a\}\st  \cc^+$
and for any $\eta \, \in\,\,D $  and  $y \in\,\,\r^n$ one has

 (4) \qquad  $|\chi_{\eta} '(y)| \le   \beta  (1 + |y|^2) e^{ - \gamma |y|^2}$,
with $\beta =\beta (a,\al) >0$, $\gamma=\gamma (a,\al) >0$.
 (positive constants only depending on $a$  and $\al $).

\noindent For $\eta \, \in\,\,D$  one has further

(5)\qquad $|\chi_{\eta}(y) -\chi_z(y)| =
|\text{\, path\,} \int _z  ^{\eta }\,\chi_z (x) d\zeta \,|  \le  |\eta  - z| \sup_{\zeta \, \in\,\,D}
 |\chi_{\zeta} ' (y)| $

 $\qquad \qquad  \le   |\eta -z| \beta (1+|y|^2) e^{- \gamma |y|^2}$.

By (2) with  $V_n : = H(\eta_n) - H(\xi_n)$,  one has to show

  (6) \qquad   $||V_n|| \to 0$   as   $n \to \infty$, with $\eta_n, \xi_n \to z$;
  where $\eta_n, \xi_n  \in D$.

\noindent  With  $h(\eta,y) : =  (\chi_{\eta}(y) - \chi_z(y))/ (\eta - z)$,  by (2) one has

  $ V_n(x) = \int_{\r^n}  (h(\eta_n,y) - h(\xi_n,y)) f(x-y) \, dy $;

   with  $v_n(y) : = h(\eta_n,y) - h(\xi_n,y)$  and  $f = wg$ with $g \in  X$, thus

  $V_n(x)  =  \int_{\r^n}  v_n(y) f(x-y)\, dy =  \int_{\r^n}  v_n(y)w(x-y)g(x-y)\,dy$,

        with   $v_n(y) \to  0$  for all $y \in \r^n$  (pointwise).

  Case $X = BUC$ or $C_0 $:
 Then $||V_n||_{wX} = ||V_n/w||_X = \sup_{x \in \r^n} || \int_{\r^n}  v_n(y)
      (w(x-y)/w(x) g(x-y) \,dy||_Y   \le
    \int_{\r^n} |v_n(y)| w(y) ||g||_X $    with (1.1).

 Now with (4) and (5) one has   $|h(\eta_n,y)| \le  \beta e^{-\gamma |y|^2}$
since $\eta_n  \in  D$, the same for $\xi_n$, so

$|v_n(y)|w(y)  \le  2 \beta e^{- \gamma |y|^2} w(y) = : h(y)$, with $h \in L^1$,
and $v_n(y) \to 0$ for each $y $; so the Lebesgue Dominated Convergence
Theorem can be applied, yielding   $V_n \to  0$  in $wX$  as desired.

 Case $X = L^p$, $1\le p<\,\infty$.

 By the above one has for all $x \in\,\,\r^n$,

$ V_n(x)/w(x) = \int_\r^n F_n(x,y) dy$,   with
 $F_n(x,y) : = v_n(y) (w(x+y)/w(x))g(x-y)$, with  $g \in\,\,L^p(\r^n,Y) $.

\noindent Here, with $||\cdot||_{p,x}  : = L^p$-norm with respect to the variable $x$,
one has  $||F_n(\cdot,y)||_{p,x} = |v_n(y)|  ||(w(x-y)/w(x))g(x-y)_{p,x} \le
  |v_n(y)|  w(y) ||g||_{L^p}$  (with $w(x-y/w(x) \le w(y)) = h(y) ||f||_{wL^p}$
 with the $h$  from above, $h \in\,\,L^1$.
 This means all the assumptions of the continuous Minkowsi inequality
 (G. Int. p.251, A 90)  are fulfilled one gets
$||V_n||_{wL^p} = ||V_n/w||_{L^p} \le \int_\r^n  |v_n(y)| w(y) ||g||_p dy
    =  ||f||_{wL^p}  \int_{\r^n}  |v_n(y)| w(y) dy $;
with the above   $ |v_n(y)| w(y) \le h(y)$ for all $y \in\,\,\r^n$ and  all  $n \in\,\N$,
 and   $h  \in\,\, L^1 $, one can apply again the Lebesgue Dominated Convergence
Theorem, getting
  $||V_n||_{wL^p}  \to 0 $ as  $n \to \,\infty$,  as desired .     \P
\enddemo

The proof is almost the same of Example 3.7.6 [A, p. 154].

Comments on your email $A_1=d/dx$

For $X = C_b(R,Y)$ , $C_{ub}$ and $AP$, with suitable $D(A)$, one has
$\sigma (A_1)  = iR$;
however, in all 3 cases I can show that $ ||R(a+it,A_1)||$ does  not tend
to $0$ as  $|t| \to \infty$, so NRE is not applicable.

Only in the case $X = P_1(R,Y) =\{f \in  C(R,Y) : f_1 = f \}$  ($f$ with period 1)
I am unable to show this, also I get only   $i 2 \pi \z  \st \sigma (A_1)  \st   i R $.

Have you any ideas about the behaviour of  $R(a+it,A_1)$  as $|t|\to \infty $

$X = P_1(R,Y) \st C_{ub}$, so the semigroup is $C_0$-semigroup and even bounded. I think the answer should follow from general theory.

The case of the Gauss-Poisson semigroup seems settled, if you checked
the proof for the case $L^p_{w_k}$.

Concerning $L^p_{w_k}$, for me the interesting case is if unbounded functions
(as in $C_{w_k}  : = w_k C_b $) are admitted, i.e.
      $ L^p_{w_k}  : = w_k  L^p : = \{f \in Y^{R^n} : f /(w_k)  \in L^p \}$.

My search for examples for NRE with $X = C_b$, $C_{ub}$, $AP$, $QP_M$, $P_a$
is now finished ($QP_M =$ quasi-periodic functions $=\{f \in AP : sp_{Bohr} f
\st M \}$, $M$ finite (or infinite)  $\st  R$ ;  $P_a : = \{f \in AP : f_a = f \}$, $0< a  \in R $,
always with $T =$ translations in $X$, with $R(t) : = R(\beta+it,A_1)$, $0\not =\beta  \in R $:

For $C_b$, $C_{ub}$, $AP$ one has    (*) inf $\{ ||R(t)|| : 0\not =t \in  R\} > 0$,
  NRE is not applicable, for no $\beta\not =0$ the $\beta+ A_1$ generates a
          holomorphic semigroup .

For $P_a$  ons has        (**)    $||R(t)||  = O(1/|t|) $   as  $|t| \to \infty$,$\beta\not =0$,
           so here all   $\beta + A_1$  generate hol. (semi)groups, NRE is
             applicable.

For $QP_M $, one has both cases :
      If $M$ contains $2$ rationally independent elements,
              (*) holds,  NRE is not applicable.

      For all other $M$,  (**) holds,  NRE is applicable .
         (But then, if  $M$ not empty and $\not = {0}$, there exists $\beta>0 $ with
            $QP_M  \st  P_{\beta} $.)

For me this was a disappointment (looking for $T$  with the $\delta$-condition
 of NRE with $\delta < 1$).
But maybe there is a  general theorem saying that if only $R(r) \to0$  as
$|t| \to \infty$, $T$ is  already holomorphic
(The best in this direction, at least in A., is Cor. 3.7.18, p. 160.)

\enddocument

Thank you for your email Poisson 10. I comment in the following:

However, A. p. 160 , Cor. 3.7.17, is indeed applicable (with e.g. a = 1),
so P is certainly holomorphic, our NRE is applicable for any $k \ge  0$.
(This would follow also already from (*) below, $z  = i t$, $t$ real., $|t|>1$.)

While working on this, I also found in A., p. 159, Cor.3.7.15;
if  I had seen this earlier, it would have saved me a lot of work and time!
However, already for k = 1 this Cor. is not applicable, the translation
group is not bounded in $C_{w_k}(R,\cc)$ .

Now the formula I mentioned above.

If  $0 \le k \in  R$, $X = UC_{w_k}(R,Y)$, $D = \{f \in C^2(R,Y) : f , f ', f '' \in  X\}$,
$A = A_2  =  (d/dx)^2$,  $z = r e^{i \phi}$,  $0< r \in R$, $|\phi| < \pi $ (I think $\pi/2$), one has :
      $\sigma (A_2)  =  C \setminus (-\infty,0]$,  and

(*) $||z R(z,A_2|| \le \beta_k (1/(\cos (\phi/2)))(1+(|z|^{1/2}\cos(\phi/2))^{-k}$
                              $\times \int^{\infty}_0 t^k e^{-t} dt$ ) ,
   with   $\beta_k =  1$ if  $0
   \le  k \le 1$, $\beta_k = 2^{k-1}$  if $k > 1$.
(So also with the new $w_k = (1 + |x|)^k $ such constants appear.)

comment:

The proof of [A, p. 154 (b)] is valid for the case $L^1_{w_k}$: The function $z\to \chi_z$ is holomorphic from $\cc^+ \to L^1_{w_k}$ if $0 < \theta < \pi/2$ for $ (\chi_z)' =d \chi_z/dz$ is continuous and belongs to $L^1_{w_k}$ for all $z\in \cc^+$ with $0 < \theta < \pi/2$. Of course $G(z)$ is not bounded in this sector.

\enddocument
$\triangledown  u \nabla u \triangle\Delta$

Problem 02.11

 Solution 1'. Take $x= at$. Then  $y= \frac{a^6 (1+t^6)}{365}$. So, if there is an integer t such that
$\frac{1+t^6}{365}= $ is also an integer, then $(ta, ma^6)$ is a pair of integers. One easily check that if $t=3$, $m=2$.

Problem 07.13.

By tangent-chord theorem and the assumption $PQ \| CA$, we get

$\angle CAB=\angle CBP= \angle BCP= \angle PQB$.

So, quadrilateral $PBQC$ is cyclic and hence $\angle PQC=\angle CBP= \angle PQB=\angle CAB$. It remains to find $\angle CAB$. Extend $CA$ to $CJ$ such that $JA= AI$. It follows $\Delta CIJ$ is congruent to $\Delta CIB$ and so $\angle AJI=\angle AIJ=(1/2)\angle ABC$. It follows $\angle CAI= 2\angle AJI=\angle ABC$ and $\angle CAB= 2\angle CAI= 2\angle ABC= 70^{\circ}$. Hence $\angle AQC= 40^{\circ}$.

\newpage

Problem 16.13

Let $T$, $S$, $D_C$ be the intersection of $m$ with  $BP$,  $AQ$ , $AB$ respectively.Let $D_R$ be the intersection of $l$ with  $AB$.  Since  $\angle PCT=\angle QCS$, $\angle CPT=\angle CQS$, we get $\angle PTC=\angle QSC$. Hence $\angle CTP=\angle STR=\angle TSR =\angle ASD_C $. Since $\angle TRS=\angle PRA$, it follows $2 \angle D_R RT=2 \angle D_CTB$. Hence $RD_R \| TD_C$.

\

\

\

\

\

\

\

\

\

\

\

\

\

\

\
\

\

\

\

\

\
\

\

\

\

\

\

02.11 (Q1), 01.13 (Q5), 02.13 (Q1-2), 05.13 (Q4-5), 07.13 (Q3), 08.13 (Q4),

 09.13(Q3), 10.13 (Q1-2), 11.13(Q2), 12.13(Q4-5), 13.13(Variant 2, Q2)

 14.13 (Q1), 15.13 (Q2-3), 16.13 (Q1).

 \

 Q1 02.11, 02.13, 10.13\qquad Q2 13.13, 11.13\qquad 07.13, 09.13

 Q4 08.13, 15.13\qquad Q5 01.13, 05.13, 12.13.

\enddocument

p. 2, line after (*) end :  instead of " Phi only e $L^{infty}(J.R)$"  write
          "and instead of Phi uniformly continuous bounded only Phi
           e $L^{infty}(J,R)$ is needed, "

p. 3. line -9, after  Phi :  replace the "and"  by  "respectively"

p. 3, line -9 , end :  replace the "and" by "respectively"

      line -8, beginning :  instead of "mollifier" write "Friedrich's mollifier"
                     ((for the reviewer))

     same line, middle :  replace the "and" by "respectively"

  line -3 :  replace the  "needed"   by "relevant"

p. 4, (1.5) ::   the "if" there seems correct

     Equation numbers (1.5), (1.4) are on p. 4 at the right side, but
        equ.numbers (1.1), (1.2) ae on the left side
      I do not know what the "correct " position should be, later
     it also varies !   (Maybe ask the editor.)

p. 7 , Ex. 2.1, line 3 :  "and so Phi  not in AP = ... " ((insert "Phi" before the
                      symbols for "not element of" ))

         Ex. 2.1, last line : the " (3.(8)"  there is correct as far as I see

p. 9, line 2 :  the bracket ")" after "...$S^1$-ap " is correct

p.11, line above Def. 3.1 :  the " b) " is correct, change the "380" to
                                   " 380/381 "
         ( the referee apparently looked only at p. 380;  if we write only 381
              instead of 380, he might say "on p. 381 there is no Prop.5.6.7.))

p. 13, 3 lines after (3.7), middle : add 3rd ")" after "..$M_kPhi$))

      3 lines above Ex. 3.3 :  ....Examples 3.3 (iii), (iv)"    ((3.3, not 3.4 ))

p. 15, 2. end : to satisfy the referee, replace " 26 pages" by " 1 - 26 "
      ((Or are the Diss Math. collected in volumes of one year , and do
       they there have different page (and volume) numbers ))

With best greetings Hans

I have another question in connection with the example of a mild
solution of    $u'' = Au + f $which is not $C^1$:
The A there $(A_2)$ is well-behaved (regular) if J = R, i.e.$ \sigma A =|=C$.
It would be even better if A is the generator of a $C_0$-semigroup.

Q : Does  $(d/dx)^2|R$ generate a $C_0$-semigroup
    (Cor.3.14.9 p. 213 [A] does not seem relevant, but I suspect
     the answer is no; do you know anything in this direction )
    Or do you have an idea how to construct a counterexample with
     A generator of a $C_0$-semigroup

Thank you for your 2 emails GEE 18, 19.   I am  now convinced that my intuition about mild solutions was completely wrong.

Concerning the example:
For $n\in \N$ , any $ X$ , $J = R$ or $R_+$,

$A_n f : = f^{(n)}$, $f\in D^n : = \{f\in C^n(J,X) : f^{(m)}\in Y \}$, $m = 0, 1,...,n$,
  with $Y : = C_{ub}(J,X)$ .

Then one has : $D^n$ linear $\st  Y$, $D^n$ dense in Y, $A_n : D^n \to Y$ linear,
                    $A_n$  is closed.

If n=1, $A_1$ generates a $C_0$-semigroup T, and so $\sigma(A_n) \not=  \cc $,
             so " A is well-behaved " .
          However the T is not holomorphic.

- Question 1 : What exactly is $\sigma(A_1)$ :

 $\sigma(A_1) =i\r$. As in  [BB, Lemma 6.3, Semigroup Forum v.54 (1997)].

 Also, this is well known result for differential operators. But in my paper I used the theorem about the boundedness of the primitives (indefinite integral) of functions.

Question 2 : Is $ \sigma(A_n) \not=  \cc$, at least for n = 2 ;
                   what exactly is $\sigma(A_2)$ :

Since $A_2= A_1^2$,  $\sigma(A_2)=\{ \la ^2\in \cc: \la\in i\r\}= -\r_+$.

Can you give me a reference where one can find results on differential operators,
especially the $A_n = (d/dx)^n$, their spectrum ,wether they
generate $C_0$-semigroups etc

Concerning the case J = compact interval, then also $D(A_n)$ and $A_n$
are defined, one easily checks that already $\sigma(A_1) = \cc$, already
$A_1$ is then "singular".

$\sigma(A_1) = \{ z \in \cc : Re z \le 0\}$  if $ J = R_+$, $A_1$  generates a $C_0$-semi-
                    group T, but $T$ is not holomorphic,$ R(z,A_1)$ does not
                    fulfil a $\delta$-conditions needed in NRE.

$\sigma(A_1)  = iR$  if J = R, $A_1$ generates even a $C_0$ group, but it is
                   not holomorphic, no "$\delta$-condition" holds .
                 (So here $\rho(A_1)$ is not a a connected set .)

Since $A_n = (A_1)^n$ (and so also for the corresponding D's), and since
      by Dunford-Schwartz I, p. 604, 10 Theorem, one has $sigma(A_1)^n$
     $ = \{z^n : z \in  \sigma(A_1)\}$, one gets e.g.

$\sigma (A_2,R_+) = \cc$ ,  this $A_2$ is singular, though it is closed.

$\sigma(A_2, R) =$ real half-line $(-\infty,0] $.

Since I do not know how the $R(.,A_2)$ for J = R behaves,  I have the
following questions :

Q1. Does $A_2$ generate a $C_0$-semigroup  (the translation-group for $A_1$
          has $A_1$ as its generator, not $A_2$.)

Q2. Does $A_2$ satisfy at least a $\delta$-condition as in NRE

\enddocument

In the following we study the following evolution differential equation

(1.1) \qquad $u^{(n)} (t)= Au (t)+\phi(t)$, $t\in J$,

\noindent where $A$ is a closed operator with dense domain $D(A) \st X$, $\rho (A)\not =\emptyset$
and $\phi\in L^{\infty}(J,X)$

We  adopt the following  definitions of mild and classical solutions.

\proclaim{Definition}1.1 $u\in C (J,X)$ is a mild solution of (1.1) if $ P^n u (t)\in D(A)$ for each $t\in J$ and there exists  $a_0, \cdots, a_{n-1} \in X$ such that $u (t)= \sum_{k=0}^{n-1} a_k t^k + P^n u (t)+ P^n \phi(t)$ for each $t\in J$.

$u\in C^{(n)} (J,X)$ is a classical solution of (1.1) if $u (t)\in D(A)$ and $u (t)$ satisfies (1.1) for each $t\in J$.
\endproclaim

\proclaim{Proposition 1.2} If  $u\in C (J,X)$ is a mild solution of (1.1) and $n > 1$, then  $  u ^{n-1} \in C(J,X)$ and $Pu(t)\in D(A)$ for each $t\in J$.
\endproclaim

\demo{Proof} Let $\la\in \rho (A)$. Then

$(A-\la)^{-1} u(t)= (A-\la)^{-1} \sum_{k=0}^{n-1} a_k t^k +(A-\la)^{-1} A P^{n}u (t) + (A-\la)^{-1}  P^{n}\phi (t)$.

\noindent Since $(A-\la)^{-1}$, $(A-\la)^{-1} A x= I+\la (A-\la)^{-1}x$, $x\in D(A)$ are bounded, then $(A-\la)^{-1} u'(t)= (A-\la)^{-1} (\sum_{k=0}^{n-1} a_k t^k )' +(A-\la)^{-1} A P^{n-1}u (t) + (A-\la)^{-1}  P^{n-1}\phi (t)$. By [1, p. 600], $(A-\la)^{-1}$  defines a one-to-one mapping of $X$ onto $D(A)$. This implies $P^{n-1} u(t)\in D(A)$, $t\in J$. It follows $ u'(t)= (\sum_{k=0}^{n-1} a_k t^k )' + A P^{n-1}u (t) +   P^{n-1}\phi (t)$. Repeating the above we conclude our statement. $\square$
\enddemo

\proclaim{Proposition 1} If $A$ is closed and $\rho (A) \not =\emptyset$, then  $A^2$ is closed  with domain $D(A^2)=\{y\in D(A): Ay \in D(A)\}$.
\endproclaim
\demo{Proof} First, we prove the case $A^{-1}$ exists. Assume that $(x_n) \st D(A^2)$ with $(x_n) \to x\in X$ and $ A^2 x_n \to z \in X$
 as $n\to \infty$. We  prove that $x\in D(A^2)$ and $z =A^2 x$.
Indeed, since  $A^{-1}$ exists and $ A^2 x_n \to z \in X$, it follows   $y_n := A^{-1} A^2x_n= A x_n   \to A^{-1}z = y \in X$.
Since $A$ is closed, it follows
$x\in D(A)$ and   $ y= Ax$. We have $ y_n \in D(A)$, $y_n \to y\in X$ and $Ay_n \to z\in X$. Again, since $A$ is closed, it follows $y\in D(A)$ and $z= Ay = A^2x$. This proves that  $A^2$ is closed with domain $D(A^2)$.

Now, assume that the inverse  $ (\la -A)^{-1}$ exists. Since $A$ is closed, $B= \la-A$ is also closed with domain $D(B)= D(A)$. By the above $B^2= \la^2 -2\la A +A^2$ is closed with domain $D(B^2)$. We prove that
$A^2 = \la^2 -2\la B + B^2 $ with domain  $D(A^2)= D(B^2)$ is closed. Assume that $(x_n) \st D(A^2)$ with $(x_n) \to x\in X$ and $ A^2 x_n \to z \in X$. We prove that $ B x_n \to y_1 \in X$ and $ B^2 x_n \to y_2 \in X$.

 \noindent Indeed, we have

 \qquad (1) $(A^2- \la^2) x_n= (-2\la B + B^2) x_n \to (z-\la^2 x )\in X$.

\noindent  It follows
$B^{-1}(-2\la B + B^2) x_n = (-2\la+ B) x_n \to B^{-1}(z-\la^2 x )\in X$. This implies

 \qquad (2) $ B x_n \to y_1= (B^{-1}(z-\la^2 x )+ 2\la x)\in X$.

\noindent  From (1) and (2) we conclude

 \qquad (3) $ B^2 x_n \to y_2= (z-\la^2 x +2\la y_1)  \in X$.

\noindent  Since $B, B^2$ are closed, we conclude $y_1= Bx$ with $x\in D(B)= D(A)$  and  $y_1 \in D(B)$ and $y_2= By_1= B^2 x$. This proves $x\in D(A^2)$ and $z= A^2 x$. Hence  $A^2$ is closed. $\square$
\enddemo

\proclaim{Proposition} If $A$ is closed and $\rho (A) \not =\emptyset$, then  $A^n$ is closed  with domain

$D(A^n)=\{y\in D(A^{n-1}): Ay \in D(A)\}$ for all $n\in \N$.
\endproclaim
\demo{Proof} First, we prove the case when $A^{-1}$ exists. Let

$(x_k) \st D(A^n)$ with $(x_k) \to x\in X$ and $ A^n x_k \to y_n \in X$
 as $k\to \infty$.

\noindent Since  $A^{-1}$ exists,  it follows for $r=1, \cdots, n-1 $ that

$y_{k,n-r }:= A^{-r} A^n x_k= A^{n-r} x_k  \in D(A^{n-r})$,  $ y_{k,n-r } \to A^{-r}y_n := y_{n-r} \in X$.

\noindent Using that  $A$ is closed for  $r= n-1, r=n-2, \cdots, r=0 $, we conclude

$x\in D(A)$,    $ y_1=  A^{-(n-1)}y_n= Ax$, $y_1\in D(A)$,  $ y_2=  A^{-(n-2)}y_n= A^2x\in D(A)$,  $\cdots$, $y_{n-1}= A^{-1} y_n = A^{n-1}x\in D(A)$ and   $ y_n= Ay_{n-1}= A^n x$. This   proves that $x\in D(A^n)$ and $y_n =A^n x$. This proves that $A^n$ is closed.

Now, assume that the inverse  $ (\la -A)^{-1}$ exists for some $\la \not =0$. Since $A$ is closed, $B= \la-A$ is also closed with domain $D(B)= D(A)$. By the above $B^n= (\la - A)^n$ is closed with domain $D(B^n)$. We prove that
$A^n = (\la - B)^n  $ is closed  with domain  $D(A^n)= D(B^n)$.

 Assume that $(x_k) \st D(A^n)$ with $(x_k) \to x\in X$ and $ A^n x_k \to z \in X$. We prove that
$ A x_k \to y_1 \in X$, $\cdots,  A^{n-1} x_k \to y_{n-1} \in X$ as $k\to \infty$.

 \noindent Indeed, we have $A^n= (\la-B)^n=\la^n+ Bp_{n-1} (B)$, where $p_{n-1} (B)=a_{n-1}+Bp_{n-2 (B)}$, $\cdots$, $p_{2} (B)=a_{2}+B p_1(B)$,  $p_{1} (B)=a_{1}+B a_0$ with $a_k \in \cc$ and  $p_{k} (B)$ are polynomials of $B$ of degree $k$.

\noindent It follows $p_{n-1} (B)x_k = B^{-1}(A^n -\la^n) x_k \to B^{-1} (z-\la^n x):= z_{n-1}\in X$, $\cdots$, $p_{1} (B)x_k = B^{-1}(p_2 (B)-a_2) x_k \to B^{-1} (z_2-a_2 x)\in X$.  This implies $ B x_k \to w_1 \in X$, $\cdots,  B^{n} x_k \to w_{n} \in X$ as $k\to \infty$. It follows $ A x_k \to y_1 \in X$, $\cdots,  A^{n-1} x_k \to y_{n-1} \in X$ as $k\to \infty$. This and the fact that $A$ is closed we conclude that $A^n$
is closed. $\square$
\enddemo

If $A =a_{i,j}$ is an  $n\times n$ square ring matrix, then (see [1, pp 58, 59])

$det \,(A)= |A|=\sum _{P} \pm a_{1,i_1}a_{2,i_2}\cdots a_{n,i_n}$, where the summation is taken over all permutations $(i_1,i_2, \cdots, i_n)$ of $(1,2, \cdots, n)$ and the sign $+$ or $-$ is taken according as the permutation is even or odd. The cofactor of the element $a_{i,j}$ is denoted by  $A_{i,j}= (-1)^{i+j}$ times the determinant of order $n-1$ that is obtained by striking out the $i^{th}$ row and $j^{th}$ column of $A$. The adjoint of $A$ is denoted by
$adj (A)= : (A_{j,i})$.
It is well known

$ a_{i,1} A_{i,1}+a_{i,2} A_{i,2}+\cdots +a_{i,n} A_{i,n}= a_{1, i} A_{1,i}+a_{2,i} A_{2,i}+\cdots +a_{n,i} A_{n,i}=|A| $,

$ a_{i,1} A_{j,1}+a_{i,2} A_{j,2}+\cdots +a_{i,n} A_{j,n}= a_{1, i} A_{1,j}+a_{2,i} A_{2,j}+\cdots +a_{n,i} A_{n,j}=0 $, $i\not = j$.

It follows

$$adj \,(A)\, A=A\, adj\, (A)= \pmatrix |A| & 0 & 0 \cdots & 0\\
0 & |A| & 0 \cdots & 0\\0 & 0 & |A| \cdots & 0\\
\vdots &\vdots &\vdots\cdots & \vdots\\ 0& 0 & 0 \cdots & 0\\ 0 &0 & 0 \cdots & |A| \endpmatrix$$.

In the following $B_n, I_n$ denote $n\times n$ matrices with

 $$B_n= \pmatrix 0 & I & 0 \cdots & 0\\
0 & 0 & I \cdots & 0\\0 & 0 & 0 \cdots & 0\\
\vdots &\vdots &\vdots\cdots & \vdots\\ 0& 0 & 0 \cdots & I\\ A &0 & 0 \cdots & 0 \endpmatrix, \qquad  I_n= \pmatrix I & 0  \cdots & 0& 0\\
0 & I  \cdots & 0 & 0\\0 & 0  \cdots & 0 & 0\\
\vdots &\vdots \cdots &\vdots & \vdots\\ 0& 0 \cdots & I & 0\\ 0 &0  \cdots & 0 & I \endpmatrix$$

Also, we denote by $M_{h}f=M_{h_{k}}\cdots M_{h_{1}}f$, $h=(h_1, \cdots, h_{k})\in \r^{k}$, $h_1 >0$, $\cdots$, $h_{k} >0$, $f\in L^1_{loc} (J,X)$.

We can adopt the definitions of mild and classical solutions of M. but without the condition uniformly continuous. Also, we can introduce  a new definition:

$u\in C^{(n-1)} (J,X)$ is a mild solution of (*) if $(u)^{(n-1)}+ x= A Pu + P\phi$.

I think we can do better than M. if we try harder.

\proclaim {Proposition} If $u\in BC(J,X)$  is  a  mild  solution  of

(*) \qquad $u^{(n)}= Au +\phi$, where

\noindent $A$ is a closed operator with dense domain $D(A) \st X$
and $\phi\in L^{\infty}(J,X)$, then for any $h \in  (\r^+)^{n+1}$ the $v = v_h  $ is  a bounded uniformly continuous solution of

(**) \qquad $v'= B_n v +F_h$, where

\noindent $v= (M_h u, (M_h u)',\cdots,(M_h u)^{(n-1)})$, $F_h = (0,\cdots,0,M_h\phi)$.

\noindent Moreover, $\sigma (B_n)=\{\la \in \cc: \la ^n\in \sigma(A)\}$.
\endproclaim

\demo{Proof} Since $u\in BC(J,X)$ is  a  mild  solution  of (*), $w= M_h u$ satisfies $w, w', \cdots, $

\noindent $w^{(n)}\in BUC (J,X)$ and $w^{(n)}= M_h u^{(n)}= M_h Au =Aw + M_h\phi$.  This shows that $v, v'\in (BUC(J,X))^n$ and $v$ is a classical solution of (**).
Here
$(BUC(J,X))^n= Y_1\times Y_2\cdots \times Y_n$ with $Y_k=BUC(J,X)$ for $k=1, \cdots, n $.

Since  det $(B_n- \la I_n)= (-\la)^n + (-1)^{n-1} A= (-1)^{n-1} (A-\la^n)$, it follows

 $\sigma (B_n)=\{\la \in \cc: \la ^n\in \sigma(A)\}$.
\enddemo

Remark. It can be shown that if $A-\la^n I$ is invertible, then the matrix, with operator elements,  $B_n -\la I_n$ has an inverse using  $(B_n -\la I_n)\cdot \text{\,\,adj\,\,} (B_n-\la I_n)= \text{\,\,adj\,\,} (B_n-\la I_n)(B_n-\la I_n)= I_n |B_n-\la I_n|$, where $|B_n-\la I_n|$ is the determinant of $B_n-\la I_n$ given by (see [1, p. 59])

$|B_n-\la I_n|= $

Example.  $$B_2-\la I_2= \pmatrix -\la I& I\\
A & -\la I \endpmatrix, \qquad (B_2-\la I_2)^{-1}= \pmatrix \frac{ -\la I}{\la^2 I-A} &\frac {-I}{\la^2 I-A}\\
\frac{-A }{\la^2 I-A}&\frac{ -\la I}{\la^2 I-A}\endpmatrix $$

\enddocument

Hello All,
Dr Sergey Ajiev, our special guest fom UNSW, would present a talk tomorrow.
The details are as follows:

Date:Wednesday, 25th Jan, 2012
Time:12:00 noon
RmNo:345, Bldg 28,Clayton

Title:

Solution to the problem of Hölder classification of infinite-dimensional spheres.

Abstract:
The uniform classification of infinite-dimensional spheres, developed
in relation with the solution of the distortion problem is more balanced
than the continuous, isometric, Lipschitz or uniform classifications of
infinite-dimensional Banach spaces. It allows to transfer a group structure,
group actions and other metric-related constructions from one space onto another.
We show that the uniformly continuous homeomorphisms can be
"upgraded" to the Hölder ones in the classical setting and establish the
explicit and, occasionally, sharp exponents of the Hölder regularity
for pairs of concrete spaces, including various Besov, Lizorkin-Triebel,
Sobolev, sequence, Schatten-von Neumann and other Banach spaces
(including lattices and more general non-commutative spaces).
Our function spaces are allowed to be anisotropic and can be
defined in terms of differences, local approximation by polynomials,
the coefficients of wavelet expansions, Littlewood-Paley
decompositions or a functional calculus. Not every equivalent
norm is geometrically friendly.
These results appear to have close ties with the presence  of a remarkable
phenomenon from the infinite-dimensional approximation theory
discovered by Tsar'kov for the uniform mappings between pairs of Lebesgue
spaces and the problems of extension and interpolation of the mappings
between the pairs of the spaces under consideration.
Among the applications of the main results are multiple examples
of spaces that do not allow any C*-algebra structure but can be endowed with
a homogeneous Hölder group structure.

\enddocument

Dear Hans,

If A is closed and $ 0\not \in   \sigma(A)$, then I can prove that $A^n$ is also closed and $0
\not \in \sigma A^n$.
In the general case, where only $\sigma A  \not = \cc$, i.e. there is $c \not \in
\sigma A$, with  $B : = c - A$,  the above shows $B^n$ is closed for all $n\in \N$,
and $D(B^n) = D(A^n)$ for all $n$ (already this was for me non-trivial).
But till now I have been unable to get $A^n$ closed.

I also can answer   your  question that if $A$ is a linear closed operator with dense domain $D(A)$, then $A^2$ Is closed with domain $D(A^2)=\{y\in D(A): Ay \in D(A)\}$.

I thought that the question is easy, but I could not answer it. I found in a book  by
\qquad  K-J Engel, R. Nagel, One -parameter semigroup for Linear Evol. Eq. , 2000  springer-Verlag  ([EN, p. 519, B.14])

\noindent  that in general $D(A^2)$ might be $0$ unless $\rho (A) \not =\emptyset$: for example $A$ generator of $C_0$-semigroup. In this case one can assume that $A^{-1}$ is invertible and prove:

\proclaim{Proposition} If $A$ is closed and $\rho (A) \not =\emptyset$, then  $A^2$ is closed  with domain $D(A^2)=\{y\in D(A): Ay \in D(A)\}$. Similarly,   $A^n$ is closed  with domain $D(A^n)=\{y\in D(A): A^{n-1}y \in D(A)\}$.
\endproclaim
\demo{Proof} Without loss of generality we can assume that $A^{-1}$ exists. Assume that $(Ax_n) \st D(A^2)$ with $y_n= Ax_n \to y \in X$ and $A y_n \to z \in X$ as $n\to \infty$. We  prove that $y\in D(A)$ and $z =Ay$.
Indeed, since  $A^{-1}$ exists, it follows   $A^{-1} A^2x_n= A x_n = y_n  \to A^{-1}z = y \in X$. Since $A$ is closed, it follows
$y\in D(A)$ and   $ z= Ay$. The case $n\in \N$ can be treated in the same way.
\enddemo

If $B= \la - A$ is closed, then $A=\la- B$ is also closed. Similarly, if $B= \la - A$, $B^2= (\la - A)^2$ are closed, then $A=\la- B$, $A^2= B^2-\la B +\la A$ are closed. the general case follows by induction.

Now, also Iam happy with the definition of mild solution of Muk. but we do not require uniform continuity. Because I can prove that for mild solutions the derivatives exist for all $k =1, \cdots , n-1$.

I hope I can implement this information soon.

Best wishes, Bolis

\enddocument
So to be able to apply this to GEE, one needs an Algebra W(to be able to consider
$\lambda - B$) which contains A and I, Ix = x for $x\in X$.
This means   W should contain all the powers $A^n$; but the $ D_n = D(A^n)$,
by Theorem 2.7 of Pazy, p. 6,
this intersection ($ = : D_0$) is even again dense in X, and furthermore one
can show that $AD_0\st D_0$ , a linear space $\st D(A) \st X$.
So with  $A_0 : = A|D_0$ and $I_0 : = I|D_0$ one has $A_0, I_0 \in \, Lin(D_0,D_0) : = \{ S : D_0 ->D_0 : S linear\}$, an algebra with unit $I_0$.

W : = smallest algebra in Lin containing $A_0$ and $I_0$ is well defined.
By the result mentioned above (resp in GEE6 below), $det\, B_0$ and  $adj B_0$
are well defined, one has  $B_0  adj B_0 =  (adj B_0) B_0 =  det (B_0)I_n$
with $B_0 : = B|D_0$  (i.e. A in B is replaced by $A_0$).
Also, $det (\lambda - B_0) = \lambda - A_0$.
I have not yet checked that an inverse for $\lambda - A_0$ gives  an
inverse for $\lambda - B_0$, and wether this gives then $\lambda$ not in
$\sigma B$  ($B_0$ is not a Banach space, just normed), and finally wether
$sigma A = sigma A_0$ (there is a relevant paper of Wrobel (a colleague
in Kiel)   which I have to study.
So there is quite some work for me before one can be certain wether this
method will work.
(If not , I have some other ideas ("direct calculation" that $ B (U adj(\lambda -B)) =
I_n$, U inverse of  $\lambda - A$, but this also would be a bit tedious).

With best greetings    Hans

In the last days I tried to settle the question of determinants, in the literature
up till now I have not found this, and p. 281 of McCloy was hard to under-
stand for me.I settled this finally with the help of the book by van der
Waerden, Moderne Algebra I, p. 57 and 58, extending and checking the
remarks there.
The result is : If R is a commutative ring (it need not contain 1, it can
have divisors of 0) then for each $n\in N$ there exist functions
                   $d_  : R^{n x n}  ->  R$
($R^{n x n}$ is the set of $n x n$  matrices with elements from the ring R)
which haave the usual properties of  determinants, especially one has
(1) p. 281  of McCloy.
If you find any place in the literature where this is formulated  a n d
shown,please let me know.
(Saying that the proof is as in the complex case is not really helpful;
the proofs I know use fields, and they are not really simple or short,
it would for me take quite some time to check everything.)

And please be patient with me, in my age progress is slow,
and by many experiences in the past I do not trust statements as
"is well known" or "similarly as in ... ", I insist on explicit proofs
or at least references.

Next I will turn to the question of what ring is suitable for the results
formulated in your GEE; I have some ideas, but again this will take
some time.

Only after this I will be able to turn to the question  what kind of "mild
solution" can be used for GEE.

\enddocument

I am not sure we can get your draft in a form that admits publication,
for me there are still a lot of open questions.

Today I was in our library to get the Jacobson book (I have not yet met
the people I wanted to asked about general determinants).
The Jacobson book is in my eyes no help, at the important points he just
says "as is well known" without giving any references - I hate such authors,
this book is no help.
In your draft you copy this, including the "It is well known" (p. 1. line 7);
I would like to know where a proof of this can be found for general
commutative rings (for fields I know several explicit references).As McCoy indicated,
for the polynomial ring $Z(x_1,...,x_m)$ one has the needed
formulas since this ring has no divisors of zero (is contained in the quotient
field of rational functions); I have not yet checke that here for the "varia-
bles" elements of an arbitrary commutative ring can be substituted - this
would yield a proof for the needed results on determinants.
Anyhow, this I think is certainly known and should be in some book,
it would save me time.

I do not agree with your " The case of operators is no problem", which
I suspect you will realize when you try to write down an explicit proof.
\enddocument

7. Proof, line 2 : I still do not see how one gets w^(n) =Aw + M_hPhi
     I only get  w^(n) = a_h  + Aw + M_h Phi with some a_h e X ,
     using the def. of mild sol.

I do not want to check the proof in M. , so I agree with you again, it is
best to drop [13].

Now some corrections for the new  tex file from today :

p. 1  4 lines above (*) :  I suggest to bring already here the corresponding
             text from the beginning  of § 2, as follows :
     ..... These examples are instructive for various conclusions concerning
         many classes of generalized almost periodic functions, they
          demonstrate ...  ((then line 3 above (*) ))
       ((The beginning of §2 can be shortened accordingly, see below.))

   (*)  could you shift the equation somewhat to the right (it looks so
          unsymmetrical to me)

   line after (*) :   ... and Phi only   $ \in L^{infty}$....     ((insert "only"  ))

 line s 2 and 3 after (*) :  instead of this (somewhat vague)  sentence I
     suggest the following :  This seems new even for uniformly continuous u.
     ((This is correct even if the result in M. should be true : in our Th.3.2
        the Phi only in $L^{infty}$, and $ J = R_+$ is admitted.))

p. 2 line 8 :  in  BAA(J,X)  and  VAA(J,X) I suggest changing the "J"  to
       " R "  ((we have not defined this; in all the other cases  the spaces are
        directly defined in the literature(not by  F|J ),  but not for BAA, VAA

p. 3, last line  :  would deleting the  " (ii) " put the remaining " ). " at the end
             of the line above  Or else delete " Phi = g and F = g AP(R,C)" ,
              just   " ..... by Ex. 3.3 (ii) ."

p. 4.  § 2 :  I suggest deleting the first "$sp_A Phi = 0$ ... ", it appears already
          in the introduction.
         If the next sentence " Ex. 2.2 is instr.... " appears in the introduction
          as I suggested, the 3rd sentence could just read  :
        " For the benefit of the reader we give first some relevant definitions:"

\enddocument
I have not yet looked at the introduction, you can send me a new version
if you want to change something.

Now some further  corrections, referring to your 124-15 of today :

p. 1, title : I suggest " On spectral criteria and inclusions  for ...."
            ((insert "crieria and" ;  only with this addition the Ex.s of §2
          appear also in the title   ))

  line above (*)  :   delete this line
        ((without it our result "appears " even more general, it holds for
         any bounded mild sol. ))

    (*)    delete  the " u(0) = x e X, "      ((see above))

p. 2, , line 6 of  § 1, end : ....  BAA(R,X)       ((the " ) " is missing

p. 3,  line 2 of (1.6) :   $ sp_F (Phi) = sp_F (Psi)$ ....    ((delete the  " : "
                   after ...(Phi) ; here this = is n o t  a definition, it is a
result
                    (which requires a proof) ))

   second line of  (1.7) : shift the "in (1.5), then ..... " to the right as in
                     (1.6), to indicate that this line belongs to (1.7)

    second line of (1.8) : shift it to the right
           ((see above; later (1.8) is used, and if a reader looks at the text as
                it is printed now he could interpret this as just the assumptions
              there are meant, overlooking the main point here, the "F c MF" ))

 p. 4  line 1 :  Begin with ´" This is "   a new line

     line 2  end :  add  "(Ex. 3.3). "    ((saving the reader the trouble of
                       trying to check this for himself))

    (1.10)   $M_hPhi = (Phi*s_h) |J$, ....       (( the " ) " is missing ))

p. 10, Ex.3.3(i), end : add after "...=R"    " , so (3.2) holds trivially. "
                   ((this now in the proof, but no one looks at the proofs))

 Proof (of Ex.3.3) (i), end :   .... = $ sp_{{0|R}} =$
       ((the subscript of should be {0|R} ,  n o t   0|R  :
          f : = 0|R  is a function but $sp_f $ is not defined;  only $sp_F $ is defined,
         where F is a subset ( "c", not " e " ) of $L^1_loc$;  here the subset
        is the linear subspace consisting of only one element , 0 ))

     line 3 of Proof of (i) :  delete the "So, (3.2) is trivially satisfied."
              ((see my remarks above))

last line end :   $= sp_{{0|R}} (g)= $    (( the "{", "}" are missing, see above))

 same line , end :   the "R " which belongs here appears on the next page,
      which is somewhat unfortunate, can you put it here ar at least
     ((with the now  longer  $ sp_rr$..... in a new line (on this page or the next)

p. 11 , References :  the " 5. Gen. EL-results..." is no longer used in this
         text as far as I see and can therefore be deleted.

\enddocument
p. 1 , title :  replace "Spctral" by "On spectral criteria and"
  ((The examples of §2 concern the "spectral criterium"  "$sp_F f$ empty
    and .... imply f e F",  therefore the "On Spectral criteria" ))

line -6 :  ... t+a e J, |f|(t) : = ||f(t)|| .    ((insert  " |f|(t) : =... "))

line -4  :  put  $1/h $  before the integral

  line after (1.5) :  begin with  "The$ sp_F(Phi)$ of .....   "  a new line

 line above (1.6) :   put a  " : "  after the (1.6)

line above (1.6) :  put the "([4, Lemmaa 1.1(C)])" at the end of the line be-
                    fore

 line after (1.7) : delete the "define", but leave the following "as" in its
         place, directly under the "If" above.

 line above (1.8) : shift the ([4, Cor.2.3(C)])" at the end of the line before

line above (1.9) :  delete the "See ([5, Prop. 2.2] , "

  same line :  put the  rest  "[6, Prop. 3.5(ii), p.421] . " at the end of the line
           before

 (1.9)  I think this (and the text following p. 4) can be deleted ;
      it will only be used in the proof of  Th. 3.3 I think, and there we just
     cite "[1,Prop. 1.1.7]"

p. 4, lines 1 -4 : see the text to (1.9) (delete this)

Prop. 1.3 :  Is this needed somewhere later, I do not remember this;
       if not, delete it ; eventually the proof can be supplied where it is
       needed, if this is only once

Prop. 1.4 could be omitted, or just a "Remark" with (i) , (ii) only,
    and without proofs .

  I will think about the ex. 1.5/1.6 .
 (I think we should make this paper as short as possible, to ease accep-
   tance.)

\enddocument

p. 11, Ex. 3.4, line 2 :  g instead of phi under the integral

    same line, end :  ... c e C, $g(t) = e^{i t^2}$.
 ((this def of g is buried somewhere in § 1, but no one will remember it))

 line -7 : after "..$g^  = R^$" insert " since $g^(s) = pi^{1/2}e^{i pi/4}
                                   e^{-i s^2 /4}$
     ((in 2009 I calculated explicitly this S'-Fourier trasnsform, I know no
      reference))

     lines -7/-6 :  delete "I follows (3.5) is trivially satisfied."
          ((this text already in Ex 3.4, end of (i) ))

   Concerning Ex 3.4 (ii), the "so (3.6) is not satisfied"  is not correct :
         (3.6) is about the homogeneous equ., then all solutions are
           $c gamma_1$, so (3.6) is true (it holds for any linear diff equ. with
          constant coefficients, from C,  general X ).

    so 3 lines above Ex.3.4, beginning, the "also" should be deleted.

   Ex. 3.4, (ii) the "$sp_0(g) = sp_F (y)$ and so (3.6) is not satisfied"
        should be replaced by "$F not in MF$"

 p. 12, line 4, end  :  ... $sp_F g  =  sp_F y $        ((y, not g))

 Question: do you have a  proof that  this F is unif. closed (follows
from the Porada inequality for $ gAP + BUC $)

 Still missing and instructive would be an example  with $Phi e F$, but
  (3.5)  false , or (3.6) false . Have you any ideas

When I looked at Ex. 3.4(ii), case (3.6), I assumed (as was the case in earlier
similar situations), that the case of the homogeneous equation was
meant, Phi = 0.
Since as I just realized here the, more general case Phi e MF is meant,
Ex. 3.4(ii) is of course correct, please disregard my corresponding
"corrections".

Additional correction :

p. 11, Ex.3.4 (ii), line 1, end :  here I would prefer $sp_F g$ (not $sp_0 g$)
    ((what has sp_0 g  to do with the questions treated here ))

       line 2 of (ii), end : here one could add " , and so F not in MF."

p. 11, line 4 :   ... Examples 3.4 (i), (ii)   ((add "s" and "(i)"

         same line    .... show....     ((delete "s")9

      next line, beginning :   in Prop. 1.1 /The.4.3 (i), (ii)...
            ((insert "Prop. 1.1"   and  "(i)" ))

    Ex. 3.4, (ii), line 1, end : ... c MF, g e F, but    ((insert "g e F, "

               same line, end :  sp_F g  =    ((instead of $sp_0 g $))

    Proof Ex.3.4, (i), line 1, end :  ...$sp_{0|R} (g)$
            ((instead of $sp_0 (g)$,  this is in 124 not defined))

      next line :   [7, (3.3)]      ((instead of (2.3) ))

    line -3, end :  .... and sp_{0|R}(g) = R .    ((instead of  sp_0 ))

      same line, end :   .... R by (i).    ((insert " by (i)) "  ))

p. 12, line 1  :   ... Theorem 4.2.1      ((not 4.1.4 ))

     line 3  replace "already"   by "in (ii) "

\enddocument
\proclaim{Lemma 2.1}  Let $f\in L^1 (\r)$  with $\widehat {f} \not =
0$ on a compact set $ K$. Then  there exists $g\in L^1 (\r)$ such that
$\widehat {g} \cdot \widehat {f} =1$ on $ K$. Moreover, one can choose $g$
such that $\widehat{g}$ has compact support and, if $f\in {\widehat\h} (\r)$,
with $g\in {\widehat\h} (\r)$.
\endproclaim

\demo{Proof} Choose a bounded open set $U$ such that $K\st U$ and
$\widehat {f} \not = 0$ on $\overline {U}$ the closure of $ {U}$. By
[16, Proposition 1.1.5 (b), p. 22], there is $k \in L^1(\r)$ such
that $\widehat {k} \cdot \widehat {f} =1$ on $ \overline {U}$. Now, choose
$\va\in \h(\r)$ such that $\va=1$ on $K$ and supp $\va\st
\overline {U}$.  Also choose $ \psi\in
\widehat{\h}(\r)$ such that $\widehat{\psi}=\va$ and take  $g= k*\psi$. Then
$\widehat {g}$ has compact support and if $f\in \widehat{\h}(\r)$, then
$\widehat{g}\in \h(\r)$ and so $g\in \widehat{\h}(\r)$. $\square$
\enddemo

\proclaim{Proposition 2.6} Let     $F\in L^{1}_{loc}  (\r,X)\cap  {\f}' (\r,X)$
then
$sp_{0,\widehat{\h}}(F)= sp_{0,\f}(F)=\text { supp }{\widehat {T}_{F}}$.
\endproclaim
\demo{Proof} Let   $\la_{0}\in \r\setminus sp_{0,\widehat{\h}}(F)$. Then there is $k\in \widehat {\h}$ such that  $\widehat {k} (\la_0)\not =0$  and $F*k=0$. We can assume that supp$\widehat {k} \st V$, where
  V  is  a compact  neighbourhood of   $\la_0$   with
$V\cap sp_{0,\widehat{\h}}(F) =\emptyset$ (otherwise we choose $\va \in \h$ such that $\va(\la_0)\not =0$ and supp$\va \st V$;  and we replace $k$ by $k*\widehat{\va}$ ). Choose  a compact neighbourhood  $W\st V$ of $\la_0$ such that $\widehat {k}\not =0 $  on $W$. By  Wiener's theorem  (=Lemma 2.1 above) there
exists    $f\in \widehat{\h} (\r)$    such  that $\hat  {f}=1$   on $W$ and $f= k*h$ with $h\in \widehat{\h} $. It follows $F*f== F*k*h=0$.
   Let   $g\in {\h}(\r)$
with  $\text {supp } g\subset  W$.   We  show   $\widehat  {T}_{F}  (g)=0$.
Indeed     $\overline{\hat  f}g=g$     so
 $ \overline {\hat{g}}^{*}=\overline {\hat{g}}^{*}*f^{*}$.  Hence
$\widehat       {T}_{F}(g)=T_{F}      (\hat  {g})=\overline{\hat
{g}}^{*}*F(0)=\overline {\hat {g}}^{*}*f^{*}*F(0)=0$,  showing   $\la_0
\not\in \text {supp  } {\widehat {T}_{F}}$.

     Conversely, suppose $\la_0 \not\in \text {supp  } {\widehat
{T}_{F}}$.  Let  $V$  be a neighbourhood of  $\la_0$    with
$V\cap \text {supp  } {\widehat {T}_{F}} =\emptyset$ and
   let $g\in \h(\r)$   with  $\text { supp  } g \subset  V$
and   $g(\la_0 )\not=0$.   Let
$f=   \overline   {\hat  {g}}\in  \widehat(\r)$    and   for   a   given    $t\in   \r$      set
$h(\la)=e^{it\la} g(\la)$.  Then
$f^{*}*F(t)=\int_{\r}\hat {g}(s-t)F(s)\,ds=T_{F}      (\hat
{h})=\widehat {T}_{F} (h)=0$.    Hence    $f*F=0$
and  $\hat  {f}(\la_0  )=\overline  {g}(\la_0  )\not=0$.    So
$\la_{0}                 \not\in                 sp_{0,\widehat{h}}(F)$.
\enddemo
\enddocument

- p. 8, 3 lines after (3.4$*$) : I do not understand this example;

by Prop. 3.3, (3.11), of [10] one has $sp_{A,S}F = sp_{A,\widehat{D}} F $: $\widehat{D}$ not $D$

for
$F \in  L^infty(R_+,X)$ and A with (3.1). So as an example for =|= one would
need either unbounded F (also interesting) or A without (3.1); $A = 0 ={0|R}$
(on your p. 3) satisfies (3.1) and cannot be used; but $ 0_+ : = {0|R_+}$
does not satisfy (3.1) and would be a candidate; however $sp_0+ $ is not
the Beurling spectrum; what $F|R_+$  can one use here

- Set $I_{L^1}(F)=\{ f\in L^1: f*F=0\}$. $sp^B (F) \st Z (\widehat{f})$ for each  $f\in I(F)$. If $f\in \h$, then $\widehat {f}$ is analytic, so $Z (\widehat{f}$ is at most countable.  Since $sp^B (F)$ is uncountable, $I(F) \cap \h = \{0\}$. Proving $sp_{0,\h} (F)=\r$ because $I_{\h} (F)=\{0\}$.

I think $\widehat {\h}$ is large enough to satisfy (3.11) of [10]. Large enough means ' $\h$ is dense in $\widehat {\h}$ or $\widehat {\h}$ is  Wiener algebra. [Reiter book: Classical Harmonic Analysis and lcg, p.22 ]. $\h$ is not a Wiener algebra.

Do you have AP88 by BP. I am trying to find it on my computer to send it to you.  If I could not find it, I will send it by mail. Also, I received 5 copies of our paper in IJPAM and I expect that they send you also 5 copies. If not I will send you 2 copies.

- p. 9, Prop. 3.1, (i), case $W = L^{\infty}$, follows from Lemma 1.1(C) of
          ABN [8].

Yes of course. But Lemma 2.3 is more interesting (and I think new) in the case $W=\f'_{ar}$.

p.18 , after (4.3) : one should mention here that the above $LF$ resp $CF$
            are holomorphic on $\cc_+$ resp $\cc\setminus iR$.  (Thank you for noticing that).

  - last line :  you mean here $L^1((\omega-\e,\omega+\e),X)$.

    Yes, as in [5, p. 324].

-p.19, line 3 : also  C((w-eps,w+eps),X).

Yes,  as in [4].

 -  (4.9)    it should be  $sp^C(F|R_+)$,  not $sp^C F$   (not defined)

 (4.9) is correct, I guess in  your question you mean $sp^L (F|R_+)$,  not $sp^L F$; we abused the notation and defined $sp^L (F): = sp^L(F|R_+)$ (in  AP125 we used

 $sp^{L^+} (F): = sp^L(F|R_+)$).

 - Prop 4.1 : For $sp^C $ (some of)  these results can be found in Pruess
        [28] p. 20, Prop. 0.4;
      I suspect that such results are also in [5] or other books, certainly
          for the Beurling spectrum.
     (Maybe one should mention this to avoid the impression that one thinks
       all this is new.)

- Yes, I think we should mention that the case $sp^C$ is not new (though Pruess is cited before).

-p. 20, 2 lines above Prop. 4.2 : the definition of $ sp_A$  seems to be
        missing in your manuscript (maybe I mentioned this already).

        (3.9)   $sp_A$, $sp_0$  seem not defined in this manuscript (resp. reference)

  $sp_0$ is defined in (1.1) and  $sp_A$ is defined 9 lines above (3.5) = 6 lines below (3.4$*$).

 -   Th. 4.3 (i), line 2 : "If also$ F \in  L^{infty}$"  I interpret as "F satisfies the
conditions of the sentence before this", and then  already $PF \in  BUC$.

    So maybe it is better to bring this in a new (i*), but without the "also"

Also here refers to the condition $0\not \in sp^L$. This is formulated by Alan and I can not challenge him in this regard. (He is on a trip in Canada  until the end of Sept.)

 - Th4.3 (ii) The conclusions imply (or are equivalent with) $F \in MAAP$
       resp. $MC_0$

- p. 21  Remark 4.4(iii) is one of the main results of this paper.
       So for the benefit of the reader could you exemplify this (at
least for me) by an explicit citation of one of the generalizations of Ingham's
Tauberian theorems  (e.g. from [5]) and then demonstrate how you get
this from your results, and even with a reduced sp instead of $ sp^{wL}$.

[5, th. 4.9.5 (Ingham)]: $F\in L^{\infty} \cap SO(R_+,X)$ and $sp^{wL} (F)=\emptyset$ implies $F\in L^1_{loc,0}(R_+,X)$ follows from Theorem 3.7 (ii) case $\A_0=C_0$.

Th.3.8(iv) : does this mean that $F \in D'_{AP}$ of our DAP, JMAA 2006

   I think the answer is no:  because $\va \in \h$ implies  $\widehat{\va}\not \in \h$ unless $\va =0$.

\enddocument

For $f: \r\to \cc $ we set $f^* (t)=\overline{f(-t)}$.

Proposition. Let $F\in L^1_{loc} (\r,X)\cap \f'(\r,X)$. Then $sp_{0,\f} (F)=$ supp $\widehat {F}$.

\demo{Proof} Let $\la_0 \in \r\setminus sp_{0,\f} (F)$. Choose  $\delta >0$ such that  $V_{\delta}=(\la_0 -\delta,\la_0 +\delta )$  satisfies $V_{\delta}\cap sp_{0,\f} (F) =\emptyset$. There is $\va\in \f(\r)$ such that
$\widehat {\va}=1 $ on $V_\delta/2$ and supp $\widehat {\va} \st V_{\delta}$. Choose $\psi\in \f(\r)$ with
supp $ {\psi} \st V_{\delta/2}$. One can show $\widehat {T}_{F} (g)=0$. Indeed, $\psi= \overline {\widehat{\va}}\cdot \psi$ and so $\overline{\widehat{\psi}}^*= \overline{\widehat{\psi}}^{*}*\va^*$. Hence
$\widehat {T}_{F} (\psi)= {T}_{F} (\widehat{\psi})=$.
\enddemo

1-p. 8, line after (3.3) : I suppose the J here is the same as in (3.1)

Yes.

2-   2 line (3.5) : in [28], Prop. 0.5 only F of polynomial growth are
         admitted, this is strictly stronger than F e AR; also, only
          $ sp^C F =$  supp$\widehat {F}$ is shown, the Beurling spectrum does
              not appear. (also not in [5] for F e AR, is then

             $ sp_{0,S}F
              =  sp^C F $ true )

I think the answer is yes even for regular temperate distribution.  I will try to send you a proof (today or tomorrow). Unfortunately, I do not know where to find  a proof in literature. But I have a proof in an old AP of BP.

- p. 9, line 6 : in [8] there is only UAA defined; do you mean here
       BAA, VAA or $BAA_u = VAA_u $
        (Most people mean BAA when they talk of aa functions.)

We  mean BAA, VAA and I think the interested reader can check that.

 -  same line :In [24] Levitan ap does not appear (they use "N-ap"),
       so just [24] does not help the reader.

       I agree, but the details require more space and time. Also reference [24] means implicitly  N-ap.

  -  Prop. 3.1(ii) :  this is a special case of our Lemma 2.2 in [8].

I do not agree in the case F in AR (R,X) which contains properly $L^{\infty}_{w_k}$.

However, I confess that I was not aware of our Lemma 2.2 in [8]. Also, the method of the proof is different.

- p. 11, Prop. 3.4(ii) :  (This is only a side remark)  This is an analog to
              Lemma 4.2 of  [10] , the case $ V = D_{L^1}$ is not included
              here, neither $AR \st  D'_{L^1}$ nor  $(D'_{L^1}) \cap L^1_{loc}  \st Ar $.

I agree that Prop 3.4 (ii) is not as general as Lemma 4.2, but can we obtain the case $J=\r_+$ of Prop 3.4 (ii) from Lemma 4.2 which proved  for distributions.

-   (3.9)   $sp_A$, $sp_0$  seem not defined in this manuscript (resp. reference)

 $sp_0$ is defined in (1.1) and  $sp_A$ is defined 9 lines above (3.5) = 6 lines below (3.4$*$).

-p. 14, end of proof for (i) :  it would help the reader if you could give a
       reference for "if all Fourier coefficients of f e EAP vanish, f e "
       (Ruess and Summers )

       I agree.

-p. 15, Th. 3.8(iii)  (side remark) Here all $ F*\psi  \in  AP$ is equivalent with
                     $F  \in  MAP$.

Yes.

     Th.3.8(iv) : does this mean that $F in D'_{AP}$ of our DAP, JMAA 2006

     Most probable. I will try to check.

-p. 17 : 2 lines above Ex.3.13 : What do you mean by "ergodicity condition"
             in Theorem 3.7(i), (ii), (v), I do not see any

The ergodicity condition of  Theorem 3.7(i), (ii), (iii) is proved at the beginning of the proof of Theorem 3.7. In part (v),  ergodicity is replaced by  the weaker condition $M_h F$ is bounded.

\enddocument

In this paper the authors study the existence of bounded  and almost automorphic solutions for the evolution equation   in a Banach space $X$ (*) $x'(t)= A x(t)
+ f(t,x(t))$ for $t\in \Bbb{R}$, where $A: D(A)\subset X\to X$ is the generator of a $C_0$-semigroup of bounded linear operators on $X$ and $f: \Bbb{R}\times X\to X$ is  an almost automorphic in $t$ uniformly with respect to the second argument. Sufficient conditions ensuring the existence of an almost automorphic solution when there is at least one bounded solution on $\Bbb{R}^+$ are given. The authors use the subvariant functional method to show that every $K$-minimizing mild solution is compact almost automorphic. Applications are provided for heat and wave equations with nonlinearities in several functional spaces.

\enddocument
In this paper is the authors study  the recurrence relation (*)\, $f_n (t)= a_n (t)f_{n-1} (t)$

\noindent $- b_n (t)\frac {\partial}{dt}f_{n-1} (t)$, for all $n\ge 1$, where $a_n $, $b_n $, $f_0 : D\to \Bbb{C}$ are analytical functions in $D$. For (*)they found an explicit formula for the sequence $f_n $ for any given $a_n $, $b_n $ $f_0 $. Also, a generalization of (*) to higher order recurrences   is discussed and a solution for special cases is given.
\enddocument
In this paper the authors  revisited some properties of pseudo-almost automorphic functions  (Theorem 2.1, Corollary 2.2, Lemma 2.3 ) which are proved for more general classes (see (2.2) in  [Basit, B., Zhang, C.,  Canad. J. Math., V. 48 (6), 1996 pp. 1138-1153], [Zhang, C., Proc. Amer. Math. Soc., 121 (1994), 167-174] and   [Zhang, C, Almost Periodic Type Functions and Ergodicity, Kluw. Acad. Publ. and Science Press, 2003 (1.5.2 pp. 65-66)]). Applications to semilinear differential equations are given in Theorem 3.2.
\enddocument

 The author continues his study of $S^p$-pseudo-almost periodicity started   in [Toka Diagana, Commun. Math. Anal. 3 (1) (2007)] (see also  mean classes in [Basit B., G\"{u}nzler, H.,  RJMP, 12(4), (2005) pp. 424-438,   introduced in (0.2)]). In this paper the existence of pseudo-almost periodic solutions  to some nonautonomous differential equations in the case when the semilinear forcing term is both continuous and  $S^p$-pseudo-almost periodic for $p >1$ is investigated. Applications to the heat equation are given.

\enddocument

 It is the aim of this paper to contribute to the structure of periodic orbits and related issues of conjugacy for the case  of endomorphisms of the $2$-torus, represented by the action (mod 1) of an integer matrix $M\in \text {Mat} (2,\Bbb {Z})$ on $\Bbb {T}^2 \simeq \Bbb{R}^2/\Bbb{Z}^2$. Counting periodic orbits on the torus  is studied, with special focus on the relation between global and local aspects and between  the dynamical zeta function on the torus and its analogue on finite lattices.The situation on the torus is determined  by the determinant, the trace and a third invariant of the matrix defining  the toral endomorphism.

 A topological group $G$ together with a compact subgroup $K$ are said to form a $Gelfand\,\,  pair$ if the set $L^1 (K\setminus G/K)$ of $K$-biinvariant integrable functions on $G$ is a commutative algebra under convolution. the authors study a family of finite Gelfand  pairs arising in connection with Heisenberg group $H=H_n (\Bbb{F})$ over finite fields of odd characteristic. The symplectic group $Sp (n,\Bbb{F})$ acts on $H$ by automorphisms. A subgroup $K$ of $Sp (n,\Bbb{F})$ yields a Gelfand  pair $(K,H)$ when the $K$-invariant functions on $H$ commute under convolution. An interesting example of this type occurs with $K$ a finite analog of the unitary group $U(n)$.
-----------------
The authors investigate the existence and uniqueness of almost periodic and almost automorphic solutions to the semilinear parabolic boundary differential equations  (SBDE)  $x'(t) = A_m x(t) + h(t,x(t))$, $t\in \Bbb{R}$, $Lx(t)= \phi (t,x(t))$, $t\in \Bbb{R}$, where $A:= A_m | \text{\,ker} L$ generates a hyperbolic analytic semigroup on a Banach space $X$ and the functions $h, \phi$ are defined on some intermediate space $X_{\beta}$, $0 < \beta < 1$ and take values in $X$ and in a boundary space $\partial X$ respectively. This paper is a continuation of the author's study  (PAMS, 134 (2006)) of  a hyperbolic differential equation, where the boundary term $\phi$ defined on the whole space $X$.

Dear Ricky,

Here is Theorem 4.3

\proclaim{Theorem 4.3} If a subset $X$ of $\Bbb{H}^{n+1}$ is compact, then it is closed and bounded.
\endproclaim
\demo{Proof} Topologically a compact set is any set for which each open cover has  a finite sub-cover. Let $\Lambda= \{N_{\epsilon}(x): x\in X\}$, where $N_{\epsilon}(x)$ is the open hyperbolic sphere  with center in $x$ and radius $\epsilon$. Then $\Lambda$ is an open cover for $X$. Let $N_{\epsilon}(x_j)$, $j=1, 2, \cdots m$ be a finite sub-cover. Since hyperbolic spheres are bounded, $X\subset \cup_{j=1}^m N_{\epsilon}(x_j) $ is bounded. Assuming that $X$ is not closed, there is $(x_k) \subset X$ which converges to $x\not\in X$. Then the open cover $\{\Bbb{H}^{n+1}\setminus \{x\}, N_{\epsilon_k}(x_k), k=1,2,\cdots\}$, where $N_{\epsilon_k}(x_k)\cap N_{\epsilon_l}(x_l)=\emptyset$, $k\not =l$ has no finite subcover. A contradiction which proves that $X$ is closed.
\enddemo
\enddocument
Let $P(z)= \sum_{j=1} ^{m} c_j\,e^{i\lambda_j}z$, $c_j\in\cc$, $\lambda_j\in\r$, $\lambda_j\not=\lambda_{j'}$ if $j\not =j'$. Hypothesis of  Lagrange on the mean motion   states that $\Delta P (t)= c\Delta t+o (\Delta t)$ when $\Delta t\to \infty$. This hypothesis is settled by
 Jessen and Tornhave (Acta Math., 77 (1945) pp. 137-279).  In this note the author simplifies this solution.
\enddocument
The purpose of this paper is to extend  some properties of Banach valued almost automorphic   functions defined on $\Bbb{R}$
to the class of almost periodic functions $AA(\Bbb{R},F)$ with $F$ $p$-Fr\'{e}chet space which is not locally convex.
Although the $p$-norm does not have all properties of an usual norm, the main properties of Banach valued almost automorphic   functions are extended to $AA(\Bbb{R},F)$. Applications to semigroups of linear operators and to dynamical systems in $p$-Fr\'{e}chet spaces are given.

\enddocument

 This paper is concerned with almost periodic functions in the sense of Besicovitch  defined on a tube $T_{K}=\{z=x+iy: x\in \r^p, y\in K\}$, where $K\st \r^p$ is the base of the tube and $p\in\N$. This work extends in some sense the results of Bauermeister (Nachr. Akad.Wiss. G\"{o}ttingen II,  Math. Phis., KI: 1972, 1975). The author found conditions for an almost periodic function in the sense of Besicovitch defined on $\r^p$ to be extended onto a tube domain as a holomorphic almost periodic function in the sense of Besicovitch.

\enddocument

 \qquad\qquad\qquad  \qquad \qquad{\bf{Report on the   paper}}

\qquad  A new approach to the spectral theory of functions and the

 \qquad\qquad\qquad  \qquad
  Loomis-Arendt-Batty-Vu theory

\qquad \qquad  \qquad \qquad \qquad\qquad\qquad by

\qquad \qquad  \qquad \qquad \qquad {\bf{ Dr Nguyen Van Minh}}

\qquad \qquad  \qquad \qquad \qquad reference number: 10.1136.

\qquad \qquad  \qquad \qquad  submitted to the Journals of the
LMS.

\

\noindent General comments.

\

\noindent This   paper is concerned with a topic involving
harmonic analysis, transform theory and operator theory.  The
results of the paper are related to Loomis's theorem [25, case
$\Bbb{X}=\cc$] on bounded uniformly continuous functions   $f:
\r\to \Bbb{X}$ and its extension to $f:\r^+\to \Bbb{X}$. See for
example [3, 6, 9, 17, 21, 24,  29, 30, 31, 32, 33, 36, 37, 40,
42].
 In the present paper the author makes some extensions of
Loomis's theorem to functions from $BC (J,\Bbb{X})$, where $J\in
\{\r^+,\r\}$.

After reading  this paper, one is left asking:

(i) Is this useful or interesting

(ii) Exactly how does this relate to the previous literature

\noindent Concerning (i), the paper offers Example 2.5 which is
already treated    (by the author and others) in [17]. But the
author did not explore Conditions $F$ and $F^+$.  For example, can
one replace item (iii) by

(iii')\qquad $f_h\in \n$ for each $f\in \n$ and $h\in \r$.

\noindent   We note that  (iii') combined with  items (i), (iv) of
condition $F$ imply  item (iii) of condition $F$. Indeed, define
for $h>0$

 $h \,M_h f (t):=\int_0^h
f(t+s)\,ds= \int_0^t [f(s+h)-f(s)]\,ds + c(h)$  with $c(h)=
\int_0^h f(s)\,ds $.

\noindent  Then $h \,M_h f- c(h)$ is a bounded primitive of
$f_h-f$, so by items (i) and (iv),  $h \,M_h f\in \n$ for each
$f\in\n$ and $h>0$. It follows $f*\phi\in \n$ for each $f\in\n$
and $\phi\in L^1(\r)$.
 Part (iii) follows because

 (*) \qquad $R(\la,\h) f = f*\phi$  for some $\phi\in
L^1(\r)$ for each $Re\,\la \not =0$.

\noindent Concerning (ii), the author is aware of the main
references for the previous works. However, the exact relation is
not obvious to the reader.

 \noindent What is the  exact relation between
 Conditions $F$   and  $F^+$   to  the notion of
 translation

 \noindent -biinvariant in [3] or $\Lambda$-classes in [5]

\noindent The proof of Lemma 2.2 is contained  in the proof of
[24, Lemma 2.5] and is a correction of the proof of [3, Lemma
2.5].

\noindent The main idea of the  proof of Lemma 2.17 is contained
in the proof of [4, Lemma 4.6.6].

\noindent  The following papers are closely related
 to this work.

[CF] Chill, R., Fasangova, E., Equality of two spectra arising in
harmonic
            analysis and semigroup theory, Proc. Amer. Math. Soc. 130 (2001),
            675-681,

 [BG] Basit, B, G\"{u}nzler, Relations between different types of spectra and
  spectral characterizations,  Semigroup Forum 76 (2008),
 217-233.

\noindent  In [BG]  the authors studied reduced spectra with
respect to classes
 $\A$ which contain all classes satisfying Conditions
 $F$ or $F^+$.

\noindent It should be noted that Corollary 2.18 of the present
paper is enough to extend  Loomis's theorem to classes satisfying
$F$ or $F^+$.

 ($\bullet$) The following example shows that Corollary
2.18 is $false$.

The function  $f:\r\to \r$, where  $f= \sin
\frac{1}{2+\cos\,t+\cos \sqrt {2} t}$ is not almost periodic but
$sp_{\n} (f)=\emptyset$ for $\n=AP(\r,\cc)$:   $f*\phi\in
AP(\r,\cc)$ for each $\phi\in L^1(\r)$ because $f*\phi\in
BUC(\r,\cc)$. In particular, by (*) above, $R(\la,D)f\in
AP(\r,\cc)$ for each $Re\, \la \not=0$.

\

\noindent Specific corrections.

\

1- p. 3, line 3: $C_0 (J,\Bbb{X}): =\{f\in
BC(J,\Bbb{X}):\lim_{t\to\infty} f(t)=0\}$ should be replaced by
$C_0 (J,\Bbb{X}):= \{f\in BC(J,\Bbb{X}):\lim_{|t|\to\infty}
f(t)=0\}$.

2. p.3, $\cc$ in (2.2), (2.3) should be replaced by $\Bbb{X}$.

3. p. 8, Lemma 2.17

In assertion (i) 'eigenvalue of $\sigma (A)$;' should be replaced
by 'eigenvalue of the operator $A$;

In line -3 '$|f(z)|$' should be replaced by '$||f(z)||$'

4. p.9, line -3, '$sp_{\n} (t)$'  should be replaced by '$sp_{\n}
(f)$'

\

 COMMENTS FOR EDITORS.

\

($\bullet\bullet$) I cannot recommend this paper for publication
in the Journals of the LMS.

\

The function $f= \sin \frac{1}{2+\cos\,t +\cos \sqrt {2} t}$ is
not almost periodic because $f$ is not uniformly continuous. But
$f$ is Stepanoff almost periodic (Levitan B. M., Almost Periodic
Functions, Moscow, 1956, p. 212) (In Russan). This implies $M_h f
\in AP$ for each $h >0$ by [7, (3.8), p. 134].

Some of the  statements of general comments are mentioned only for
the editors:  The argument (iii') implies (iii) and the argument
that $f= \sin \frac{1}{2+\cos\,t+\cos \sqrt {2} t}$ is not almost
periodic but has void almost periodic spectrum. I think that the
author would identify my identity if he reads these statements.

 Bolis Basit

 e-mail address: bolis.basit\@sci.monash.edu.au

Date of completion of report: August, 11, 2008.

\enddocument

1.  How strong is the case for publication of the paper weak.

2.  If the paper has serious defects, how strong would the case be
if a successful revision were carried out  weak. END OF QUESTIONS
On completion, please send this report to publications\@lms.ac.uk.

\enddocument
--------------------------------------------------------------------------------------------------

EXPLANATORY NOTES AND GUIDELINES

The Editorial Board is anxious to reach decisions as quickly as
possible. Firstly, please acknowledge receipt of this form by
e-mail to publications@lms.ac.uk and let me know your estimate of
when you will return the report to me. The default time, when we
send out the first reminder, varies according to the journal:
Bulletin = 6 weeks; Journal = 8 weeks and Proceedings = 10 weeks,
but we can change this according to your estimate. Please send
your report by e-mail to publications@lms.ac.uk.

 You should receive thanks and an acknowledgement from us within seven days of the receipt of your report.
If you do not hear from us, the report may have been lost in
transit, so please then contact me.

 Here is some more general guidance on what we are looking for:

 A referee should assess the originality, correctness, importance and interest of an article.
Ultimately, responsibility for correctness lies with the author
but we do ask that referees evaluate the correctness of the major
results. Editors also expect to receive views on matters such as
the balance between achievement and length, how a paper compares
with other work on related subjects published in leading journals,
and the quality of exposition.  Referees should bear in mind that
the Society receives significantly more good material than it has
resources to print.  In view of the competition for publication
space, the Editors must turn away some well-recommended items.
Objective and carefully reasoned referees' reports are therefore
greatly valued by the  Editors. The LMS will respect the
confidentiality of your identity as a referee and you should not,
under any circumstances, communicate directly with the author(s)
yourself. Under data protection legislation, any part of the
report may be requested by the authors but information will only
be given after all possible features have been removed from the
report that could identify you as the referee. Parts of the report
that are considered helpful to the authors in preparing a better
paper may be sent to them by the Editors.

 If, for any reason, you are now unwilling or unable to review the paper, please let us know immediately.
In this case, we would be grateful if you could suggest another
suitable referee.

 Yours sincerely,

 Ben Holmes
 Assistant Editor

 On completion, please send this report to publications@lms.ac.uk.

\enddocument

In this paper the authors studied two species time-delayed
predator-prey lotka-Voltera type dispersal systems with periodic
coefficients in which the prey species can disperse among $n$
patches, while the density-independent predator species is
confined to one of patches cannot disperse (*)
$\frac{dx_1(t)}{dt}= x_1(t)[a_1 (t)-
b_1(t)x_1(t)-c(t)\int_{-\infty}^{0} k_{1,2}(s)y(t+s)\, ds] +$
$\sum_{j=1}^n [\alpha_{1,j}(t)d_{1,j}(t) x_j
(t-\tau_{1,j})-d_{j,1}(t) x_1(t)]$,\,\,
 $\frac{dx_i(t)}{dt}= x_i(t)[a_1 (t)- b_i(t)x_i(t)]+
\sum_{j=1}^n [\alpha_{i,j}(t)d_{i,j}(t) x_j
(t-\tau_{i,j})-d_{j,i}(t) x_i(t)]$, $ i= 2,\cdots, n$,\,\,

\noindent $\frac{dy(t)}{dt}=y(t)[- e(t)-f(t)\int_{-\infty}^{0}
k_{2,1}(s)x_1(t+s)\, ds] $, $t\in [0,\infty)$, $x_i$ denote the
population density of prey species in the $i$ patch and $y$ is the
population density of predator species. Sufficient conditions on
the boundedness, permanence and existence of positive periodic
solution for (*) is established. The theoretical results is
confirmed by a special example and numerical simulation.

\enddocument
Let $\sc^d$ be the unit sphere in $\r^{d+1}$ and $L^p (\sc^d)$ the
Banach space of measurable functions  $f: \sc^d\to \cc$ such that
$||f||_p= (\int_{\sc^d} |f(z)|^p d\sigma(z))^{1/p}<\infty$ if
$1\le p<\infty$ and $||f||_{\infty}=$ sup$_{z\in \sc^d} |f(z)|
<\infty$. The author found necessary density conditions for
Marcinkiewicz-Zygmond inequalities and interpolation  for the
space of spherical harmonics in $\sc^d$ with respect to the $L^p$
norm. It is proved that there are no complete interpolation
families for $p\not =2$.

\
\
\

The goal of this paper is to investigate the asymptotic behaviour
of systems governed by infinite delay differential equations in
terms of the attractors of associated truncated finite delay
equations and their numerical approximations. The authors proved
the upper semicontinuous convergence of approximate attractors for
an infinite delay equation  of logistic type, first for the
associated truncated delay equation with finite delay and then for
numerical scheme applied to the truncated equation.

\
Let $X$ be a commutative complex $C^*$-algebra with unit $e$ and
$\Bbb {T}$, $D$ denote the unit circle respectively the open unit
disk of $\cc$. Set $L^{\infty}(\Bbb {T},X)=\{f: \T\to X$
measurable, $||f||_{\infty}= \text{ess\,\,sup}_{0\le \theta\le
\pi} ||f(e^{i\theta})|| <\infty \}$. In this note the author
discussed the topological properties of the maximal ideal space of
the Banach algebra $L^{\infty}(\Bbb {T},X)$. A property for the
extreme points of the closed unit ball in the Banach space
$H^{\infty}(\Bbb {T},X)$ is given.

Let $D$ be a compact subset of $\r^n$,  $\r^{k\times m}$ be the
space of $k\times m$ matrices, $AP (\r^{k\times m})$ the space of
almost periodic functions $h: \r\to \r^{k\times m}$ and $AP (D,
\r^{k\times m})$ the space of almost periodic functions $g:
\r\times D\to \r^{k\times m}$. If $h\in AP (D,\r^{k\times m})$,
Mod$h$ will denote the minimal subgroup of $\widehat {\r}$
containing all the characters of the fourier series of $h$. The
author studied the linear system (*) $\dot{x} =A(t)x$, $x\in D$,
$t\in \r$, where $A\in \r^{n\times n}$.
 Let $M_1= $ Mod$\{aij: i\le j, 1\le i,j \le n\}$, $M_2= $ Mod$\{aij: i > j, 1\le i,j \le n\}$
 and Mod$(A)= M_1\oplus M_2$. The existence problem of irregular
 with respect to $M_1$  almost periodic solutions are studied. If
 (*) has such a nontrivial solution, then some of the diagonal
 elements of $A(t)$ are constant.
\enddocument
The aim of this paper is to investigate the convergence of the
solutions of the following system of delay differential equations
(1.1)\,\,\,$x'_1 (t)= -F_1 (x_1(t))+F_1 (x_2(t-r_2))$, $x'_2 (t)=
-F_2 (x_2(t))+F_2 (x_1(t-r_21))$, where $r_1$ and $r_2$ are given
positive constants, $F_1$ and $F_2$ are continuous and strictly
increasing on $\r$. It is shown that each solution of (1.1) tends
to a constant vector as $t\to \infty$.  This  extends the existing
ones in the literature.

\

\

In this paper the authors studied permanence of dynamical systems
which of interest in mathematical ecology.  Sufficient conditions
are given for the permanence and global attractivity of positive
equilibria of multi-species competition system with delay. The
results of the paper improve and extend some earlier
investigations of other authors.

\

\

Let $G$ be a locally compact group and $H$  its closed subgroup.
 Let $P(G)$ and  $P(H)$ be respectively the sets of continuous positive definite
functions on $G$ and $H$.
 $H$   is called extending if every function $ \phi\in P(H)$
 extends to some function in $P(G)$. $G$ has the extension
 property if each closed subgroup of $G$ is extending. Denote by
 $P_H(G)=\{\phi\in P(G): \phi=1 $ on $H\}$. $H$ is separating if
 for each $x\in G\setminus H $ there is $\phi \in P_H(G)$ such
 that $\phi(x)\not = 1$. In this paper the authors continued to investigate  earlier results
on  related aspects of these extending  and separating
properties.
\enddocument

Let $S$ be a locally compact topological semigroup and $M_a (S)$
be the space of all measures $\mu\in M(S)$(the space of bounded
Radon measures) for which the mappings $x\mapsto |\mu|* \delta_x$
and $x\mapsto \delta_x *|\mu|$ from $S\to M(S)$ are weakly
continuous. If $\cup \{\text{supp} \mu: \mu\in M_a (S)\}$ is dense
in $S$, then $S$ is called foundation semigroup. The author
studies the question of whether a closed subspace of $LUC(S)^*$ or
(respectively $M_a (S)$) invariant under $L_x^{**}$ (respectively
$L_x$), $x\in S$ ($S$ a foundation semigroup), is a left ideal.

In this paper the authors discuss almost periodicity of solutions
for evolution equation in Banach spaces (*) $u'(t)+ A(t) u(t)=
f(t)$, where $A(t): X\to X$ is an operator for $t\in \r$, $f:
\r\to X$ is a continuous function and $X$ is a Banach space. The
main problem here is: Under what conditions is any bounded
solution of (*) almost periodic. For special Hilbert space
(compact synchronous), some results on almost periodicity of all
solutions are established. In particular, the results of Haraux on
$\r^2$ [J. Diff. Eq., 66 (1987), 62-70] are extended to $\r^n$.

deal with the asymptotic behavior, the nonoscillation and the
stability of solutions of linear neutral delay difference
equations with periodic coefficients and constant delays (E)
$\Delta (x_n + \sum_{i \in I} \, c_i x_{n-\sigma_i}) = a(n) x_n +
\sum_{j\in J} \, b_j x_{n-\tau_j} $. The results are achieved via
a positive root of an associated equation which is, in a sense,
the corresponding characteristic equation. Let $G_{\r}$ be a real
semisimple Lie group, possibly not compact. One of the most
important examples of equivariantly closed forms is the symplectic
volume form $ d\beta$ of a coadjoint orbit $\Omega$. Even if
$\Omega$ is not compact, the integral $ \int_{\Omega}\, d\beta$
exists as a distribution on the Lie algebra $\frak{g}_{\r}$. This
distribution is called the Fourier transform of the coadjoint
orbit. In this paper, the author applies his localization results
(J. Funct. Anal. 203 (2003) no.1, 197-236, 215 (2004)) no.1, 50-66
and Progr. Math., 220, MA, 2004) to get a geometric derivation of
Harish-Chandra's formula for the Fourier transform of regular
semisimple coadjoint orbits. An explicit computation is made for
the coadjoint orbits of elements of $\frak{g}^*_{\r}$ which are
dual to regular semisimple elements lying in a maximally split
Cartan subalgebra of $\frak{g}_{\r}$.

Let $G$  be not simply connected  group and  $L(G)$ be the loop
group consisting of connected components. Let $\widehat{L(G)}$
denote a central extension and $ \widetilde{\cc}^*$ be a finite
covering of ${\cc}^*$ acting on $\widehat{L(G)}$ covering the
natural action of ${\cc}^*$-action on $L(G)$. The goal of this
paper is to give an explicit formula for the character $\chi_V$
restricted to the connected components of $\widehat{L(G)}\rtimes
\widetilde{\cc}^*$ which do not contain the identity element. This
gives a generalization of the Kac-Weyl character formula.

In this paper the authors study the existence of almost
automorphic mild solutions to the perturbed abstract differential
equation (*)$\dot {x}(t)= (A+B)x(t)$, where $A$ is a generator of
$C_0$-group and $B:X\to X$ a bounded operator. However, the
condition (4) $AA(X)\cap L^1(\Bbb{R}, R(B))$ is dense in
$CB(\Bbb{R}, R(B))$ is possible if and only if $B=0$.

In this paper the authors established a maximal characterization
for Stepanoff bounded functions $B S^p (\r)$ analogous  to a well
known  maximal characterization for $L^p (\r)$ when $1< p< infty$.
However, the  method of the resembles  those used for results on
maximal functions.

The authors used I for denoting  1 and also for denoting an
interval of $\r$.
\bigskip

\bigskip

\bigskip

In this paper the authors continue some earlier  studies of
Helgason-Fourier transform of $L^1$ and $L^p$ functions on
Riemannian symmetric spaces X. Analogues of the Hausdorff-Young
and Hardy-Littlewood inequalities, the Wiener Tauberian theorem
and some uncertainty theorems on Riemannian symmetric spaces of
noncompact type are formulated.

\enddocument